\documentclass{article}
\usepackage[cp1251]{inputenc}
\usepackage[english]{babel}
\usepackage{amsmath}
\usepackage{amsthm}
\usepackage{mathtext}
\usepackage{amssymb}
\usepackage{amscd}
\usepackage{graphics}
\usepackage{euscript}
\usepackage{bbm}
\usepackage[matrix,arrow,curve]{xy}
\theoremstyle{definition}
\newtheorem{defin}{Definition}
\newtheorem*{example}{Example}
\newtheorem*{rem}{Remark}
\newtheorem{quest}{Question}
\theoremstyle{plain}
\newtheorem{thm}{Theorem}
\newtheorem*{thm*}{Theorem}
\newtheorem{prop}{Proposition}
\newtheorem{cor}[prop]{Corollary}
\newtheorem*{prop*}{Proposition}
\newtheorem{lemma}{Lemma}
\DeclareMathOperator{\Qsym}{Qsym} \DeclareMathOperator{\Sym}{Sym}
 \DeclareMathOperator{\id}{id}
\DeclareMathOperator{\Ker}{Ker}\DeclareMathOperator{\im}{Im}
 
\DeclareMathOperator{\link}{link} \DeclareMathOperator{\Hom}{Hom}
\DeclareMathOperator{\Ls}{Ls}
 \DeclareMathOperator{\Coh}{H}
\DeclareMathOperator{\C}{H} \DeclareMathOperator{\pt}{pt}
\DeclareMathOperator{\Lyn}{LYN}

\newcommand {\ib}[1]{\textit{\textbf{#1}}}

\begin{document}
\title{Ring of Polytopes, Quasi-symmetric functions\\ and Fibonacci numbers}
\author{Victor M.~Buchstaber \and Nickolai Erokhovets.}
\date{}
\maketitle
\begin{abstract}
In this paper we study the ring $\mathcal{P}$ of combinatorial convex polytopes.
We introduce the algebra of operators $\mathcal{D}$ generated by the
operators $d_k$ that send an $n$-dimensional polytope $P^n$ to the sum
of all its $(n-k)$-dimensional faces. It turns out that $\mathcal{D}$
is isomorphic to the universal Leibnitz-Hopf algebra with the antipode
$\chi(d_k)=(-1)^kd_k$. Using the operators $d_k$ we build the generalized
$f$-polynomial, which is a ring homomorphism from $\mathcal{P}$ to the
ring $\Qsym[t_1,t_2,\dots][\alpha]$ of quasi-symmetric functions with coefficients in $\mathbb Z[\alpha]$.
The images of two polytopes coincide if and only if their flag $f$-vectors are equal.
We describe the image of this homomorphism over the integers and prove that over
the rationals it is a free polynomial algebra with dimension of the $n$-th
graded component equal to the $n$-th  Fibonacci number. This gives a representation
of the Fibonacci series as an infinite product. The homomorphism is an isomorphism on the graded group
$BB$ generated by the polytopes introduced by Bayer and Billera to find the linear
span of flag $f$-vectors of  convex polytopes. This gives the group $BB$ a structure of the ring isomorphic to $f(\mathcal{P})$.
We show that the ring of polytopes has a natural Hopf comodule structure over the Rota-Hopf algebra of posets.
As a corollary we build a ring homomorphism $l_{\alpha}\colon\mathcal{P}\to\mathcal{R}[\alpha]$ such that $F(l_{\alpha}(P))=f(P)^*$,
where $F$ is the Ehrenborg quasi-symmetric function.
\end{abstract}
\section{Introduction}
Convex polytopes is a classical object of convex geometry.
In recent times the solution of problems of convex geometry involves results of algebraic geometry and topology, commutative and homological algebra. There are remarkable results lying on the crossroads of the polytope theory, theory of complex manifolds, equivariant topology and singularity theory (see the survey \cite{BR}). The example of successful collaboration of the polytope theory and differential equations can be found in \cite{Buch}.

To study the combinatorics of convex polytopes we develop the approach based on the ring $\mathcal{P}$ of combinatorial convex polytopes. In the focus of our interest we put the ring of linear operators $\mathcal{L}(\mathcal{P})=\Hom_{\mathbb Z}(\mathcal{P},\mathcal{P})$, and its subring generated by the operators $d_k,\,k\geqslant 1$ that map an $n$-dimensional polytope to the sum of all it's $(n-k)$-dimensional faces. $d=d_1$ is a derivation, so $\mathcal{P}$ is a differential ring.

The ring $\mathcal{D}$ has a canonical structure of the Leibnitz-Hopf algebra (see Section 4). We prove that $\mathcal{D}$ is isomorphic to the universal Hopf algebra in the category of Leibnitz-Hopf algebras with the antipode $\chi(d_k)=(-1)^kd_k$.
$$
\mathcal{D}\simeq\mathcal{Z}/J_{\mathcal{D}}
$$
where $\mathcal{Z}=\mathbb Z\langle{Z_1,Z_2,\dots\rangle}$ is the universal Leibnitz-Hopf algebra -- a free associative Hopf algebra with the comultiplication
$$
\Delta Z_k=\sum\limits_{i+j=k}Z_i\otimes Z_j,
$$
and $J_{\mathcal{D}}$ is a two-sided Hopf ideal generated by the relations $\sum\limits_{i+j=n}(-1)^iZ_iZ_j=0,\quad n\geqslant 1$.

This algebra appears in various application of theory of Hopf algebras in combinatorics: over the rationals it is isomorphic to the graded dual of the odd subalgebra $S_{-}(\Qsym[t_1,t_2,\dots],\zeta_{\mathcal{Q}})$ \cite[Remark 6.7]{ABS}, to the algebra of forms on chain operators \cite[Proposition 3.2]{BLiu}, to the factor algebra of the algebra of Piery operators on a Eulerian poset by the ideal generated by the Euler relations \cite[Example 5.3]{BMSW}.

The action of $\mathcal{D}$ on $\mathcal{P}$ satisfies the property $D_{\omega}(PQ)=\mu\circ(\Delta D_{\omega})(P\otimes Q)$ so the ring $\mathcal{P}$ has the structure of a Milnor module over $\mathcal{D}$.

The graded dual Hopf algebra to $\mathcal{Z}$ is the algebra of quasi-symmetric functions $\Qsym[t_1,t_2,\dots]$.
The $\alpha$-character $\xi_{\alpha}\colon\mathcal{P}\to\mathbb Z[\alpha]$ induces the homomorphism $\mathcal{P}\to\mathcal{D}^*[\alpha]$, and the composition with the inclusion $\mathcal{D}^*\subset\Qsym[t_1,t_2,\dots]$ gives the ring homomorphism
$$
f:\mathcal{P}\to\Qsym[t_1,t_2,\dots][\alpha]
$$
We call $f(\alpha,t_1,t_2,\dots)$ the generalized $f$-polynomial. In can be shown that $f(\alpha,t_1,0,0,\dots)=f_1(\alpha,t)$ is a homogeneous $f$-polynomial in two variables introduced in \cite{Buch}. \\For simple polytopes $f(\alpha,t_1,t_2,\dots)=f_1(\alpha,t_1+t_2+\dots)$.

This paper has the following structure:

Sections 2-4 are devoted to the necessary definitions and important facts about quasi-symmetric functions and Hopf algebras.

Section 5 contains topological realizations of the Hopf algebras we study.

Section 6 recalls important facts about flag $f$-vectors and emphasizes some peculiarities we need in this paper.

Section 7 is devoted  to the cone and the bipyramide operators $C$ and $B$ used by M.~Bayer and L.~Billera \cite{BB} to find the linear span of flag $f$-vectors of polytopes. We show that the operations $C$ and $B$ can be defined on the ring $\Qsym[t_1,t_2,\dots][\alpha]$ in such a way that $f(CP)=Cf(P)$ and $f(BP)=Bf(P)$.

In Section 8 the structure theorem for the ring $\mathcal{D}$ is proved.

In Section $9$ we introduce the generalized $f$-polynomial and find the functional equations that describe the image of the ring $\mathcal{P}$ in $\Qsym[t_1,t_2,\dots][\alpha]$. These equations are equivalent to the generalized Dehn-Sommerville equations discovered by M.~Bayer and L.~Billera in \cite{BB}.

In Section $10$ the condition that determines the generalized $f$-polynomial in a unique way is found.

In Section $11$ we describe the graded dual Hopf algebra to the Hopf algebra $\mathcal{D}$ as a subalgebra of $\Qsym[t_1,t_2,\dots]$ and find the image of $\mathcal{P}$ in the ring $\mathcal{D}^*[\alpha]$.

In Section $12$  the multiplicative structure of the ring $f(\mathcal{P}\otimes \mathbb Q)$ is described.  We prove that $f(\mathcal{P}\otimes \mathbb Q)$ is a free polynomials algebra with dimension of $n$-th graded component equal to the $n$-th Fibonacci number $c_n$ ($c_0=c_1=1$, $c_{n+1}=c_n+c_{n-1}$). This gives the decomposition of the Fibonacci series into the infinite product
$$
\frac{1}{1-t-t^2}=\sum\limits_{n=0}^{\infty}c_nt^n=\prod\limits_{i=1}^{\infty}\frac{1}{(1-t^i)^{k_i}},
$$
where $k_i$ is the number of generators of degree $i$. The infinite product converges absolutely in the interval $|t|<\frac{\sqrt{5}-1}{2}$.
The numbers $k_n$ satisfy the inequalities $k_{n+1}\geqslant k_n\geqslant N_n-2$, where $N_n$ is the number of the decompositions of $n$ into the sum of odd numbers.

In Section $13$ we introduce the multiplicative structure on the graded group $BB$ generated by the Bayer-Billera polytopes arising from the isomorphism with the ring $f(\mathcal{P})$.

Nowadays Hopf algebras is one of the central tools in combinatorics.  There is a well-known Hopf algebra of posets introduces by Joni and Rota in \cite{JR}. Various aspects of this algebra were studied in \cite{Ehr,ABS,Sch1,Sch2}. A generalization of the Rota-Hopf algebra was proposed in  \cite{RS}.
In \cite{Ehr} R.~Ehrenborg introduced the $F$-quasi-symmetric function, which gives a Hopf algebra homomorphism from the Rota-Hopf algebra to $\Qsym[t_1,t_2,\dots]$.
In Section $14$ we show that there is a natural right Hopf comodule structure over the Rota-Hopf algebra on the ring $\mathcal{P}$. It gives the natural ring homomorphism $l_{\alpha}:\mathcal{P}\to\mathcal{R}[\alpha]$. We prove that $F(l_{\alpha}(P^n))=f(P^n)^*$ and show that the Hopf comodule structure agree with the homomorphisms $f^*$ and $F$.

Section $15$ deals with the problem of description of flag $f$-vectors of polytopes. For simple polytope  $g$-theorem gives the full description of the set of flag $f$-vectors. As it is mentioned in \cite{Z1} even for $4$-dimensional non-simple polytopes the corresponding problem is extremely hard. We
show how this problem can be stated in terms of the ring of  polytopes.

\section{The Ring of Polytopes}
For the notion of a convex polytope see the books \cite{Gb,Z2}.
\begin{defin}
A {\itshape ring of polytopes} $\mathcal{P}$ is a ring of all combinatorial convex
polytopes. As an abelian group\linebreak
$\mathcal{P}=\mathbb Z+\bigoplus\limits_{n\geqslant 1,m-n\geqslant 1}
P^{n,\,m-n}$, where $\mathcal{P}^{n,\,m-n}$ is a free abelian group generated by all $n$-dimensional
combinatorial polytopes with $m$ facets. The direct product of
polytopes $P\times Q$ defines a multiplication and the point
$1=\pt$ is a unit.
\end{defin}

It can be proved that $\mathcal{P}$ is a ring of polynomials
generated by all indecomposable polytopes. $n$ and $m-n$ are
graduations of the ring and it easy to see that there
is a finite number of $n$-dimensional polytopes with $m$ facets,
so an abelian group $\mathcal{P}^{n,\,m-n}$ is finitely generated.

\begin{defin}
Let us denote $\mathcal{P}^0=\mathbb Z$, $\mathcal{P}^n=\bigoplus\limits_{m=n+1}^{\infty}\mathcal{P}^{n,\,m-n},\, n\geqslant 1$,  and $\mathcal{P}^{[n]}=\bigoplus\limits_{k=0}^{n}\mathcal{P}^n$.

We will use the notations: $P$ for a polytope, $P^n$ for a polytope of dimension $n$, $p$ for an element of $\mathcal{P}$, $p^n$ for an element of $\mathcal{P}^n$.
\end{defin}

\begin{defin}
Let us define the ring of linear operators  $\mathcal{L}(\mathcal{P})=\Hom_{\mathbb Z}(\mathcal{P},\mathcal{P})$
with the multiplication given by a composition of operators.
\end{defin}

\begin{defin}
Let $d_k$ be an operator defined on each polytope $P^n$ as the sum of
all $(n-k)$-dimensional faces.
$$
d_kP^n=\sum\limits_{F^{n-k}\subset P^n}F^{n-k}
$$
For example, $d_0=\id$, and $d=d_1$ is a facets operator. Since
$d(P\times Q)=(dP)\times Q+P\times (dQ)$, it is a derivation. So
the ring of polytopes is a differential ring.

Let us denote by $\mathcal{D}=\mathcal{D}(\mathcal{P})\subset\mathcal{L}(\mathcal{P})$ the subring
of operators on $\mathcal{P}$ generated by $d_k$, $k\geqslant 1$.
\end{defin}

The ring of polytopes contains a differential subring
$\mathcal{P}_s$, generated by all simple polytopes.

In the ring of simple polytopes we have the relations:
$\left.d_k\right|_{\mathcal{P}_s}=\left.\frac{d^k}{k!}\right|_{\mathcal{P}_s}$,
so $\mathcal{D}(\mathcal{P}_s)$ is isomorphic to the divided power
ring $d_id_j={i+j\choose i}d_{i+j}$.

In the general case as an abelian group the ring $\mathcal{D}$ is
generated by the monomials $D_{\omega}=d_{j_1}\dots d_{j_k}$ where
$\omega$ is a $k$-tuple $(j_1,\dots,j_k)$, for some $k\geqslant
0,\;j_i\geqslant 1$.

\begin{defin}
A {\itshape composition} $\omega$ of a number $n$ is an ordered
set $\omega=(j_1,\dots,j_k),\;j_i\geqslant 1$, such that
$n=j_1+\dots+j_k$. Let us denote $|\omega|=n$.
\end{defin}
The number of compositions of $n$ into exactly $k$ parts is given
by the binomial coefficient ${n-1\choose k-1}$.

It is easy to see that
$$
d_k(P\times Q)=\sum\limits_{i+j=k}d_iP\times d_jQ=d_kP\times \
Q+d_{k-1}P\times dQ+\dots+P\times d_kQ.
$$
and
\begin{equation}\label{DD}
D_{\omega}(P\times
Q)=\sum\limits_{\Omega'+\Omega''=\omega}
D_{\omega'}P\times D_{\omega''}Q.
\end{equation}
where for the compositions $\omega=(j_1,\dots,j_k)$, $\omega'=(j_1',\dots,j_{l'}')$, $\omega''=(j_1'',\dots,j_{l''}'')$, $\Omega'$ and $\Omega''$ are all the $k$-tuples such that
$$
\Omega'=(0,\dots,j_1',\dots,0,\dots,j_{l'}'\dots,0),\quad \Omega''=(0,\dots,j_1'',\dots,0,\dots,j_{l''}''\dots,0),
$$

For example,
\begin{multline*}
dd_3(P\times
Q)=(dd_3P)\times(Q)+(dd_2P)\times(dQ)+(d^2P)\times(d_2Q)+(dP)\times
(d_3Q)+\\
+(d_3P)\times(dQ)+(d_2P)\times(d^2Q)+(dP)\times(dd_2Q)+P\times(dd_3Q);
\end{multline*}
The list of summands corresponds to the decompositions:
\begin{multline*}
(1,3)=(1,3)+(0,0)=(1,2)+(0,1)=(1,1)+(0,2)=(1,0)+(0,3)=\\
=(0,3)+(1,0)=(0,2)+(1,1)=(0,1)+(1,2)=(0,0)+(1,3).
\end{multline*}
This formula gives a coproduct in $\mathcal{D}$  defined on the generators by the formula
$$
\Delta d_k=\sum\limits_{i+j=k}d_i\otimes d_j
$$
and thus $\mathcal{D}$ obtains the structure of a graded Hopf algebra with $\deg d_k=k$.

\begin{defin}[\cite{N}]
A \emph{Milnor module} $M$ over the Hopf algebra $X$ is an algebra
with unit $1\in k$ which is also a module over $X$ satisfying
\[ x(uv) = \sum x'_n(u)x''_n(v), \quad x\in X,\; u,v \in M, \;\, \Delta x = \sum x'_n \otimes x''_n. \]
\end{defin}

By definition of the comultiplication on the ring $\mathcal{D}$
$$
\mu\Delta(D_{\omega})(P\otimes Q)=D_{\omega}(PQ),
$$
therefore for the ring $k=\mathbb Z$ we obtain the proposition:

\begin{prop}
The ring of polytopes $\mathcal{P}$ is a Milnor module over the Hopf algebra $\mathcal{D}$.
\end{prop}

Another example of operators in $\mathcal{L}(\mathcal{P})$ gives the multiplication by the elements of $\mathcal{P}$:
$$
[P](Q)=PQ
$$
\begin{prop}
The following relation holds
$$
d_k[P]=\sum\limits_{i=0}^k[d_iP]d_{k-i}
$$
\end{prop}
\begin{proof}
$$
d_k(PQ)=\sum_{i=0}^{k}(d_iP)d_{k-i}Q
$$
\end{proof}
Thus any operator in the ring $\mathcal{P}\mathcal{D}$ generated by $\mathcal{D}$ and $\{[P],P\in\mathcal{P}\}$ can be expressed as a sum of operators
$$
\sum\limits_{P,\,\omega}[P]D_{\omega}
$$

\begin{defin}
Define an {\itshape $\alpha$-character}  $\xi_{\alpha}:\mathcal{P}\to\mathbb
Z[\alpha]$ by the formula
$$
\xi_{\alpha}P^n=\alpha^n.
$$

For $\alpha=0$ set
$$
\xi_0(P^n)=
\begin{cases}
1,&P^n=\pt;\\
0,&n\geqslant 1.
\end{cases}
$$
Let us denote $\varepsilon=\xi_0$.
\end{defin}
\begin{defin}
Let us define an operator $\Phi:\mathcal{P}\to\mathcal{P}[t]$,
$\Phi(t)\in\mathcal{D}[[t]]$:
$$
\Phi(t)=1+dt+d_2t^2+d_3t^3+\dots=\sum\limits_{k=0}^{\infty}d_kt^k.
$$
\end{defin}

\begin{prop}
$$
\xi_{-\alpha}\Phi(\alpha)=\xi_{\alpha}\label{Euler}
$$\label{Euler}
\end{prop}
\begin{proof}
Let us note that $\xi_{\alpha}d_kP^n=f_{n-k}(P^n)$ -- the number of $(n-k)$-dimensional faces.
Then for any polytope $P^n$ we have
$$
\xi_{-\alpha}\Phi(\alpha)P^n=(-\alpha)^n+(-\alpha)^{n-1}f_{n-1}\alpha+\dots+f_0\alpha^n=\left(\sum\limits_{i=0}^n(-1)^if_i\right)\alpha^n=\alpha^n=\xi_{\alpha}P^n,
$$
since the Euler formula gives the relation:
$$
f_0-f_1+\dots+(-1)^{n-1}f_{n-1}+(-1)^n=1
$$
\end{proof}
\begin{defin}
Let $P^n$ be the $n$-dimensional polytope
$$
P^n=\{\ib{x}\in\mathbb R^n:\langle \ib{a}_i,\ib{x}\rangle+1\geqslant
0\}
$$
Then a {\itshape dual {\upshape (or {\itshape polar})} polytope}
is defined as $P^*=\{\ib{y}\in\mathbb R^n:\langle
\ib{y},\ib{x}\rangle+1\geqslant 0,\;\forall \ib{x}\in P^n\}$.

There is a bijection $F\longleftrightarrow F^{\diamondsuit}$
between the $i$-faces of $P$ and the $(n-1-i)$-faces of $P^*$ such
that
$$
F\subset G\Leftrightarrow G^{\diamondsuit}\subset
F^{\diamondsuit}.
$$
The operation $*$ defines a linear operator on the ring
$\mathcal{P}$ that does not belong to $\mathcal{P}\mathcal{D}$
\end{defin}

\section{Quasi-Symmetric Functions}

\begin{defin} Let $t_1,t_2,\dots$ be a finite or an infinite set of
variables. For a composition $\omega=(j_1,\dots,j_k)$ consider a
{\itshape quasi-symmetric monomial}
$$
M_{\omega}=\sum\limits_{l_1<\dots<l_k}t_{l_1}^{j_1}\dots
t_{l_k}^{j_k}.
$$
\end{defin}
Degree of the monomial $M_{\omega}$ is equal to
$|\omega|=j_1+\dots+j_k$.

For any two monomials $M_{\omega'}$ and $M_{\omega''}$ their product in the ring of polynomials $\mathbb Z[t_1,t_2,\dots]$ is equal to
$$
M_{\omega'}M_{\omega''}=\sum\limits_{\omega}\left(\sum\limits_{\Omega'+\Omega''=\omega}\right)M_{\omega},
$$
where for the compositions $\omega=(j_1,\dots,j_k)$, $\omega'=(j_1',\dots,j_{l'}')$, $\omega''=(j_1'',\dots,j_{l''}'')$, $\Omega'$ and $\Omega''$ are all the $k$-tuples such that
$$
\Omega'=(0,\dots,j_1',\dots,0,\dots,j_{l'}'\dots,0),\quad \Omega''=(0,\dots,j_1'',\dots,0,\dots,j_{l''}''\dots,0),
$$

This multiplication rule of compositions is called {\itshape the
overlapping shuffle multiplication}.

For example,
\begin{enumerate}
\item
$$
M_{(1)}M_{(1)}=\left(\sum\limits_{i}t_i\right)\left(\sum\limits_{j}t_j\right)=\sum\limits_{i}t_i^2+2\sum\limits_{i<j}t_it_j=M_{(2)}+2M_{(1,\,1)}.
$$
This corresponds to the decompositions
$$
(2)=(1)+(1),\; (1,1)=(1,0)+(0,1)=(0,1)+(1,0).
$$
\item
\begin{multline*}
M_{(1)}M_{(1,\,1)}=\left(\sum\limits_{i}t_i\right)\left(\sum\limits_{j<k}t_jt_k\right)=\sum\limits_{i<j}t_i^2t_j+\sum\limits_{i<j}t_it_j^2+3\sum\limits_{i<j<k}t_it_jt_k=\\
=M_{(2,\,1)}+M_{(1,\,2)}+3M_{(1,\,1,\,1)}.
\end{multline*}
This corresponds to the decompositions
\begin{gather*}
(2,1)=(1,0)+(1,1),\;(1,2)=(1,0)+(1,1),\\
(1,1,1)=(1,0,0)+(0,1,1)=(0,1,0)+(1,0,1)=(0,0,1)+(1,1,0).
\end{gather*}
\item
\begin{multline*}
M_{(1,\,1)}M_{(1,\,1)}=\left(\sum\limits_{i<j}t_it_j\right)\left(\sum\limits_{k<l}t_kt_l\right)=\sum\limits_{i<j}t_i^2t_j^2+2\sum\limits_{i<j<k}t_i^2t_jt_k+\\
+2\sum\limits_{i<j<k}t_it_j^2t_k+2\sum\limits_{i<j<k}t_it_jt_k^2+6\sum\limits_{i<j<k<l}t_it_jt_kt_l=\\
=M_{(2,\,2)}+2M_{(2,\,1,\,1)}+2M_{(1,\,2,\,1)}+2M_{(1,\,1,\,2)}+6M_{(1,\,1,\,1,\,1)}.
\end{multline*}
This corresponds to the decompositions
\begin{gather*}
(2,2)=(1,1)+(1,1),\;(2,1,1)=(1,1,0)+(1,0,1)=(1,0,1)+(1,1,0),\\
(1,2,1)=(1,1,0)+(0,1,1)=(0,1,1)+(1,1,0),\;(1,1,2)=(1,0,1)+(0,1,1)=(0,1,1)+(1,0,1),\\
(1,1,1,1)=(1,1,0,0)+(0,0,1,1)=(1,0,1,0)+(0,1,0,1)=(1,0,0,1)+(0,1,1,0)=\\
=(0,1,1,0)+(1,0,0,1)=(0,1,0,1)+(1,0,1,0)=(0,0,1,1)+(1,1,0,0).
\end{gather*}
\end{enumerate}

Thus finite integer combinations of quasi-symmetric monomials form
a ring. This ring is called {\itshape a ring of quasi-symmetric
functions} and is denoted by $\Qsym[t_1,\dots,t_n]$ (where $n$
is the number of variables. In the case of an infinite number of
variables it is denoted by $\Qsym[t_1,t_2,\dots]$).

\begin{prop}
A polynomial $g\in\mathbb Z[t_1,\dots,t_m]$ is a finite linear
combination of quasi-symmetric monomials if and only if
$$
g(0,t_1,t_2,\dots,t_{m-1})=g(t_1,0,t_2,\dots,t_{m-1})=\dots=g(t_1,\dots,t_{m-1},0)
$$
\end{prop}
\begin{proof}
For the quasi-symmetric monomial $M_{\omega}$ we have
$$
M_{\omega}(t_1,\dots,t_i,0,t_{i+1},\dots,t_{m-1})=
\begin{cases}
0,&\omega=(j_1,\dots,j_m);\\
M_{\omega}(t_1,\dots,t_i,t_{i+1},\dots,t_{m-1}),&\omega=(j_1,\dots,j_k),\;k<m.
\end{cases}
$$
So this property is true for all quasi-symmetric functions.

On the other hand, let the condition of the proposition be true.
Let us prove that for a fixed composition $\omega$ any two
monomials $t_{l_1}^{j_1}\dots t_{l_k}^{j_k}$ and
$t_{l_1'}^{j_1}\dots t_{l_k'}^{j_k}$ have the equal coefficients
$g_{l_1,\,\dots,\,l_k}^{j_1,\,\dots,\,j_k}$ and
$g_{l_1',\,\dots,\,l_l'}^{j_1,\,\dots,\,j_k}$.

Let $l_i+1<l_{i+1}$ or $i=k$ and $l_i<m$. Consider the
corresponding coefficients in the polynomial equation:
$$
g(t_1,\dots,t_{l_i-1},0,t_{l_i},t_{l_i+1},\dots,t_{m-1})=g(t_1,\dots,t_{l_i-1},t_{l_i},0,t_{l_i+1},\dots,t_{m-1}).
$$
On the left the monomial $t_{l_1}^{j_1}\dots
t_{l_i}^{j_i}t_{l_{i+1}-1}^{j_{i+1}}\dots t_{l_k-1}^{j_k}$ has the
coefficient
$g_{l_1,\,\dots,\,l_i+1,\,l_{i+1},\dots,\,l_k}^{j_1,\,\dots,\,j_k}$
and on the right
$g_{l_1,\,\dots,\,l_i,\,l_{i+1},\,\dots,\,l_k}^{j_1,\dots,j_k}$, so
they are equal. Now we can move the index $l_i$ to the right to
$l_i+1$ and in the same manner we can move $l_i$ to the left, if
$l_{i-1}<l_i-1$ or $i=1$ and $l_i>1$.

Now let us move step by step the index $l_1$ to $1$, then the
index $l_2$ to $2$, and so on. At last we obtain that
$g_{l_1,\,\dots,\,l_k}^{j_1,\,\dots,\,j_k}=g_{1,\,\dots,\,k}^{j_1,\,\dots,\,j_k}$.
\end{proof}
The same is true in the case of an infinite number of variables,
if we consider all the expressions
$$
r+\sum\limits_{k=1}^{\infty}\sum\limits_{1\leqslant
l_1<\dots<l_k}\sum\limits_{j_1,\,\dots,\,j_k\geqslant 1}r_{
l_1,\,\dots,\,l_k}^{j_1,\,\dots,\,j_k}t_{l_1}^{j_1}\dots
t_{l_k}^{j_k}.
$$
of bounded degree: $|\omega|=j_1+\dots+j_k<N$. As a corollary we
obtain another proof of the fact that quasi-symmetric functions
form a ring.

In \cite{H} M.~Hazewinkel proved the Ditters conjecture that
$\Qsym[t_1,t_2,\dots]$ is a free commutative algebra of
polynomials over the integers.

The $n$-th graded component $\Qsym^n[t_1,t_2,\dots]$ has rank
$2^{n-1}$.

The numbers $\beta_i$ of the multiplicative generators of degree $i$ can be found by a recursive relation:
$$
\frac{1-t}{1-2t}=\prod\limits_{i=1}^{\infty}\frac{1}{(1-t^i)^{\beta_i}}
$$

There are ring homomorphisms:

\begin{gather*}
v_{m+1}:\Qsym[t_1,\dots,t_m]\to\Qsym[t_1,\dots,t_{m+1}]:\quad
v_{m+1}M_{\omega}(t_1,\dots,t_m)=M_{\omega}(t_1,\dots,t_{m+1})\\
u_m:\Qsym[t_1,\dots,t_{m+1}]\to\Qsym[t_1,\dots,t_m]:\quad
t_{m+1}\to 0;
\end{gather*}
The mapping $u_m$ sends all the monomials $M_{\omega}$
corresponding to the compositions $(j_1,\dots,j_{m+1})$ of length
$m+1$ to zero, and for the compositions $\omega$ of smaller length
$u_mM_{\omega}(t_1,\dots,t_{m+1})=M_{\omega}(t_1,\dots,t_m)$.

It is easy to see that $u_mv_{m+1}$ is the identity map of the
ring $\Qsym[t_1,\dots,t_m]$.

Given $m>0$ there is a projection $\Pi_{\Qsym}:\mathbb
Z[t_1,\dots,t_m]\to\Qsym[t_1,\dots,t_m]\otimes\mathbb Q$:\\
for $j_1<j_2<\dots<j_k$
$$
\Pi_{\Qsym} t_{j_1}^{j_1}\dots
t_{j_k}^{j_k}=\frac{1}{{m\choose k}}\left(\sum\limits_{l_1<\dots<l_k}t_{l_1}^{j_1}\dots
t_{l_k}^{j_k}\right),
$$
which gives the average value over all the monomials of each type.

Then for the quasi-symmetric monomial
$M_{\omega},\;\omega=(j_1,\dots,j_k)$ we have:
$$
\Pi_{\Qsym} M_{\omega}=\sum\limits_{1\leqslant j_1<\dots<j_k\leqslant
m}\frac{1}{{m\choose k}}M_{\omega}={m\choose k}\frac{1}{{m\choose k}}M_{\omega}=M_{\omega},
$$
So $\Pi_{\Qsym}$ is indeed a projection.
\begin{rem} We see that in the theory of symmetric and quasi-symmetric
functions there is an important additional graduation -- the
number of variables in the polynomial.
\end{rem}

\begin{defin}
For a composition $\omega=(j_1,\dots,j_k)$ let us define the composition $\omega^*=(j_k,\dots,j_1)$.
\end{defin}
The correspondence $M_{\omega}\to (M_{\omega})^*=M_{\omega^*}$ defines an involutory ring homomorphism
$$
*:\Qsym[t_1,t_2,\dots]\to \Qsym[t_1,t_2,\dots].
$$

\section{Hopf Algebras}
In this part we follow mainly the notations of \cite{H}. See also \cite{BR} and \cite{CFL}.

Let $R$ be a commutative associative ring with unity.
\subsection{Leibnitz-Hopf Algebras}
\begin{defin}
A {\itshape Leibnitz-Hopf algebra} over the ring $R$ is an
associative Hopf algebra $\mathcal{H}$ over the ring $R$  with a
fixed sequence of a finite or countable number of multiplicative
generators $H_i$, $i=1,2\dots$ satisfying the comultiplication
formula
$$
\Delta H_n=\sum\limits_{i+j=n}H_i\otimes H_j,\quad H_0=1.
$$

A {\itshape universal Leibnitz-Hopf algebra} $\mathcal{A}$ over
the ring $R$ is a Leibnitz-Hopf algebra with the universal
property: for any Leibnitz-Hopf algebra $\mathcal{H}$ over the
ring $R$ there exists a unique Hopf algebra homomorphism
$\mathcal{A}\to\mathcal{H}$.

Consider the free associative Leibnitz-Hopf algebra over the
integers $\mathcal{Z}=\mathbb Z \langle Z_1,Z_2,\dots\rangle$ in
countably many generators $Z_i$.
\begin{prop}
$\mathcal{Z}$ is a universal Leibnitz-Hopf algebra over the
integers, and $\mathcal{Z}\otimes R$ is a universal Leibnitz-Hopf
algebra over the ring $R$.
\end{prop}

Set $\deg Z_i=i$. Let us denote by $\mathcal{M}$ the graded dual Hopf
algebra over the integers.
\end{defin}

It is not difficult to see that $\mathcal M$ is precisely the
algebra of quasi-symmetric functions over the integers. Indeed,
for any composition $\omega=(j_1,\dots,j_k)$ we can define
$m_{\omega}$ by the dual basis formula
$$
\langle m_{\omega},Z^{\sigma}\rangle=\delta_{\omega,\,\sigma},
$$
where $Z^{\sigma}=Z_{r_1}\dots Z_{r_l}$ for a composition
$\sigma=(r_1,\dots,r_l)$. Then the elements $m_{\omega}$ multiply
exactly as the quasi-symmetric monomials $M_{\omega}$.

Let us denote
$$
\Phi(t)=1+Z_1t+Z_2t^2+\dots=\sum\limits_{k=0}^{\infty}Z_kt^k
$$
Then the comultiplication formula is equivalent to
$$
\Delta\Phi(t)=\Phi(t)\otimes\Phi(t)
$$

Let $\chi:\mathcal{Z}\to\mathcal{Z}$ be the {\itshape antipode},
that is a linear operator, satisfying the property
$$
1\star\chi=\mu\circ1\otimes\chi\circ\Delta=\eta\circ\varepsilon=\mu\circ\chi\otimes
1\circ\Delta=\chi\star 1,
$$
where $\varepsilon:\mathcal{Z}\to\mathbb Z$, $\varepsilon(1)=1$,
$\varepsilon(Z^{\sigma})=0$, $\sigma\ne\varnothing$, is a
{\itshape counit}, $\eta:\mathbb Z\to\mathcal{Z}$, $\eta(a)=a\cdot 1$ is a {\itshape unite} map, and $\mu$ is a
multiplication in $\mathcal{Z}$.

Then $\Phi(t)\chi(\Phi(t))=1=\chi(\Phi(t))\Phi(t)$ and
$\{\chi(Z_n)\}$ satisfy the recurrent formulas
\begin{equation}\label{X1}
\chi(Z_1)=-Z_1;\quad
\chi(Z_{n+1})+\chi(Z_n)Z_1+\dots+\chi(Z_1)Z_n+Z_{n+1}=0,\;n\geqslant
1.
\end{equation}

In fact, since $\Phi(t)=1+Z_1t+\dots$, we can use the previous
formula to obtain
\begin{gather*}
\chi(\Phi(t))=\frac{1}{\Phi(t)}=\sum\limits_{i=0}^{\infty}(-1)^i\left(\Phi(t)-1\right)^i;\\
\chi(Z_n)=\sum\limits_{m=1}^n(-1)^m\sum\limits_{a_1+\dots+a_m=n}d_{a_1}\dots
d_{a_m}\quad a_i\geqslant 1,
\end{gather*}
This formula defines $\chi$ on generators, and therefore on the whole algebra.

The homomorphism $R:\mathcal{Z}\to\mathcal{D}:Z_k\to
d_k$, defines in the ring of polytopes the structure of a right module over the Leibnitz-Hops algebra $\mathcal{Z}$.
Let us mention, that the word $Z^{\omega}=Z_{j_1}\dots Z_{j_k}$ corresponds to the operator $D_{\omega^*}=d_{j_k}d_{j_{k-1}}\dots d_{j_1}$.

\begin{defin}
A {\itshape universal commutative Leibnitz-Hopf algebra}
$\mathcal{C}=\mathbb Z[C_1,C_2,\dots]$ is a free commutative polynomial Leibnitz-Hopf algebra in
generators $C_i$ of degree $i$. We have $\mathcal{C}=\mathcal{Z}/J_{\mathcal{C}}$, where the ideal $J_{\mathcal{C}}$ is generated by the relations
$Z_iZ_j-Z_jZ_i$.

It is a self-dual Hopf algebra and the graded dual Hopf algebra is naturally
isomorphic to the algebra of symmetric functions $\mathbb
Z[\sigma_1,\sigma_2,\dots]=\Sym[t_1,t_2,\dots]\subset\Qsym[t_1,t_2,\dots]$
generated by the symmetric monomials
$$
\sigma_i=M_{\omega_i}=\sum\limits_{l_1,\,\dots,\,l_i}t_{l_1}\dots
t_{l_i},
$$
where $\omega_i=\underbrace{(1,\;\dots,\;1)}_{i}$.

The isomorphism $\mathcal{C}\simeq \mathcal{C}^*$ is given by the
correspondence $C_i\to \sigma_i$.
\end{defin}

\subsection{Lie-Hopf Algebras}
\begin{defin}
A {\itshape Lie-Hopf algebra} over the ring $R$ is an associative
Hopf algebra $\mathcal{L}$ with a fixed sequence of a finite or
countable number of multiplicative generators $L_i$, $i=1,2\dots$
satisfying the comultiplication formula
$$
\Delta L_i=1\otimes L_i+L_i\otimes 1,\quad L_0=1.
$$
A {\itshape universal Lie-Hopf algebra} over the ring $R$ is a
Lie-Hopf algebra $\mathcal{A}$ satisfying the universal property:
for any Lie-Hopf algebra $\mathcal{H}$ over the ring $R$ there
exists a unique Hopf algebra homomorphism
$\mathcal{A}\to\mathcal{H}$.

Consider the free associative Lie-Hopf algebra over the integers
$\mathbb Z\langle U_1,U_2,\dots\rangle$ in countably many
variables $U_i$.
\begin{prop}
$\mathcal{U}$ is a universal Lie-Hopf algebra, and
$\mathcal{U}\otimes R$ is a universal Lie-Hopf algebra over the
ring $R$.
\end{prop}

Let us set $\deg U_i=i$. Then the graded dual Hopf algebra is denoted
by $\mathcal{N}$. This is the so-called {\itshape shuffle
algebra}.
\end{defin}

The antipode $\chi$ in $\mathcal{U}$ has a very simple form
$\chi(U_i)=-U_i$.

\subsection{Lyndon Words}
A well-known theorem in the theory of free Lie algebras (see, for example, \cite{Re}) states
that the algebra $\mathcal{N}\otimes\mathbb Q$ is a commutative
free polynomial algebra in the so-called {\itshape Lyndon words}.
\begin{defin}
Let us denote by $[a_1,\dots,a_n]$ an element of $\mathbb N^*$, that is the word over $\mathbb N$, consisting of symbols $a_1,\dots,a_n$, $a_i\in\mathbb N$. Let us order the words from $\mathbb N^*$ lexicographically, where any symbol is larger than nothing,
that is $[a_1,\dots,a_n]>[b_1,\dots,b_m]$ if and only if there is an $i$ such that $a_1=b_1,\dots,a_{i-1}=b_{i-1}, a_i>b_i, 1\leqslant i\leqslant\min\{m,n\}$, or $n>m$ and $a_1=b_1,\dots,a_m=b_m$.

A {\itshape proper tail} of a word $[a_1,\dots,a_n]$ is a word of the form $[a_i,\dots,a_n]$ with $1<i\leqslant n$. (The empty word and one-symbol word have no proper tails.)

A word is {\itshape Lyndon} if all its proper tails are larger than the word itself. For example, the words $[1,1,2]$, $[1,2,1,2,2]$, $[1,3,1,5]$ are Lyndon and the words $[1,1,1,1]$, $[1,2,1,2]$, $[2,1]$ are not Lyndon. The set of Lyndon words is denoted by $\Lyn$.

The same definitions make sense for any totally ordered set, for example, for the set $\{1,2\}$ or for the set of all odd positive integers.
\end{defin}

The role of Lyndon words is described by the following theorem:
\begin{thm*}[Chen-Fox-Lyndon Factorization \cite{CFL}] Every word $w$ in $\mathbb N^*$ factors uniquely into a decreasing concatenation product of Lyndon words
$$
w=v_1*v_2*\dots*v_k,\quad v_i\in\Lyn,\;v_1\geqslant v_2\geqslant\dots\geqslant v_k
$$
\end{thm*}
For example, $[1,1,1,1]=[1]*[1]*[1]*[1]$, $[1,2,1,2]=[1,2]*[1,2]$, $[2,1]=[2]*[1]$.

The algebra $\mathcal{N}$ is additively generated by the words in $\mathbb N^*$. The word $w=[a_1,\dots,a_n]$ corresponds to the function
$$
\langle [a_1,\dots,a_n], U^{\sigma}\rangle=\delta_{w,\,\sigma}
$$
where $\sigma$ is a composition $\sigma=(r_1,\dots,r_l)$, $U^{\sigma}=U_{r_1}\dots U_{r_l}$, and for the word $[a_1,\dots,a_n]$
$$
\delta_{w,\,\sigma}=\begin{cases}
1,&\sigma=(a_1,\dots,a_n);\\
0,&\mbox{else}
\end{cases}
$$
The multiplication in the algebra $\mathcal{N}$ is the so-called {\itshape shuffle multiplication}:
$$
[a_1,\dots,a_n]\times_{sh}[a_{n+1},\dots,a_{m+n}]=\sum\limits_{\sigma}[a_{\sigma(1)},\dots,a_{\sigma(n)},a_{\sigma(n+1)},\dots,a_{\sigma(n+m)}],
$$
where $\sigma$ runs over all the substitutions $\sigma\in S_{m+n}$ such that
$$
\sigma^{-1}(1)<\dots<\sigma^{-1}(n)\mbox{ and }\sigma^{-1}(n+1)<\dots<\sigma^{-1}(n+m).
$$
For example,
\begin{gather*}
[1]\times_{sh}[1]=[1,1]+[1,1]=2[1,1];\\
[1]\times_{sh}[2,3]=[1,2,3]+[2,1,3]+[2,3,1];\\
[1,2]\times_{sh}[1,2]=[1,2,1,2]+[1,1,2,2]+[1,1,2,2]+[1,1,2,2]+[1,1,2,2]+[1,2,1,2]=2[1,2,1,2]+4[1,1,2,2].
\end{gather*}
There is a well-known shuffle algebra structure theorem:
\begin{thm*}
$\mathcal{N}\otimes \mathbb Q=\mathbb Q[\Lyn]$, the free commutative algebra over $\mathbb Q$ in the symbols from $\Lyn$.
\end{thm*}

The proof follows from the following theorem concerning shuffle products in connection with Chen-Fox-Lyndon factorization.
\begin{thm*}
Let $w$ be a word on the natural numbers and let $w=v_1*v_2*\dots *v_m$ be it's Chen-Fox-Lyndon factorization. Then all words that occur with nonzero coefficient in the shuffle product $v_1\times_{sh}v_2\times_{sh}\dots\times_{sh}v_m$ are lexicographically smaller or equal to $w$, and $w$ ocurs with a nonzero coefficient in this product.
\end{thm*}

Given this result it is easy to prove the shuffle algebra structure theorem in the case of arbitrary totally ordered subset $M=\{m_1,m_2,\dots\}\subset\mathbb N$. Let us denote by $\mathcal{U}_{M}$ the free associative Lie-Hopf algebra $\mathbb Z\langle U_{m_1},U_{m_2},\dots\rangle$, let $\mathcal{N}_{M}$ be its graded dual algebra, and $\Lyn_M$ be the corresponding set of all Lyndon words. We need to prove that $\mathcal{N}_M\otimes\mathbb Q=\mathbb Q[\Lyn_M]$.

Let $m_1$ be the minimal number in $M$. The smallest word $[m_1]$ is Lyndon.

Given a word $w$ we can assume by induction that all the words lexicographically smaller than $w$ have been written as polynomials in the elements of $\Lyn_M$. Take the Chen-Fox-Lyndon factorization $w=v_1*v_2*\dots*v_m$ of $w$ and consider, using the preceding theorem, $$
v_1\times_{sh}v_2\times_{sh}\dots\times_{sh}v_m=aw+(\mbox{reminder}).
$$
By the theorem the coefficient $a$ is nonzero and all the words in (reminder) are lexicographically smaller than $w$, hence they belong to $\mathbb Q[\Lyn_M]$. Therefore $w\in\mathbb Q[\Lyn_M]$. This proves generation. Since each monomial in $n$-th graded component of $\mathcal{N}\otimes \mathbb Q$ has a unique Chen-Fox-Lyndon decomposition, the number of monomials in $\mathbb Q[\Lyn_M]$ of graduation $n$ is equal to the number of monomials in $\mathcal{N}_{M}\otimes\mathbb Q$, monomials in Lyndon words are linearly independent. This proves that Lyndon words are algebraically independent.

\begin{cor}
The graded dual algebras to the free associative Lie-Hopf algebras $\mathcal{U}_{12}=\mathbb Z\langle U_1,U_2\rangle$ and $\mathcal{U}_{odd}=\mathbb Z\langle U_1,U_3,U_5,\dots\rangle$ are free polynomial algebras in Lyndon words $\Lyn_{12}$ and $\Lyn_{odd}$ respectively.\label{Poly}
\end{cor}
The correspondence
$$
1+Z_1t+Z_2t^2+\dots+=\exp(U_1t+U_2t^2+\dots);\\
$$
defines an isomorphism of Hopf algebras $\mathcal{Z}\otimes\mathbb
Q\simeq\mathcal{U}\otimes\mathbb Q$.

However it is not true that $\mathcal{N}$ is a free
polynomial commutative algebra over the integers.

\section{Topological realization of Hopf algebras}
In \cite{BR} A.~Baker and B.~Richter showed that the ring of
quasi-symmetric functions has a very nice topological
interpretation.

We will consider $CW$-complexes $X$ and their homology $\C_*(X)$
and cohomology $\C^*(X)$ with integer coefficients.

Let us assume that homology groups $\C_*(X)$ have no torsion.

The diagonal map $X \to X \times X$ defines in $\C_*(X)$ the structure of a graded coalgebra with the comultiplication
$\Delta \colon \C_*(X) \to \C_*(X)\otimes \C_*(X)$ and the dual structure of a graded algebra in $\C^*(X)$ with the multiplication
$\Delta^* \colon \C^*(X)\otimes \C^*(X) \to \C^*(X)$.

In the case when $X$ is an $H$-space with the multiplication $\mu\colon X\times X\to X$ we obtain the Pontryagin product $\mu_*\colon \C_*(X)\otimes \C_*(X)\to\C_*(X)$ in the coalgebra $\C_*(X)$ and the corresponding structure of a graded Hopf algebra on $\C_*(X)$. The cohomology ring $\C^*(X)$ obtains the structure of a graded dual Hopf algebra with the diagonal mapping $\mu^* \colon \C^*(X) \to \C^*(X)\otimes \C^*(X)$.

For any space $Y$ the loop space $X=\Omega Y$ is an $H$-space. A continuous mapping $f\colon Y_1\to Y_2$ induces a mapping of $H$-spaces $\Omega f\colon X_1 \to X_2$, where $X_i=\Omega Y_i$. Thus for any space $Y$ such that $X=\Omega Y$ has no torsion in integral homology we obtain the Hopf algebra $\C_*(X)$. This correspondence is functorial, that is any continuous mapping $f\colon Y_1\to Y_2$ induces a Hopf algebra homomorphism $f_* \colon \C_*(X_1) \to \C_*(X_2) $

By the Bott-Samelson theorem \cite{BS}, $\C_*(\Omega\Sigma X)$ is
the free associative algebra $T(\widetilde\C_*(X))$ generated by $\widetilde \C_*(X)$.
This construction is functorial, that is a continuous mapping $f:X_1\to X_2$
induces a ring homomorphism of the corresponding tensor algebras arising from
the mapping $f_*\colon\widetilde \C_*(X_1)\to\widetilde\C_*(X_2)$.
Denote elements of $\C_*(\Omega\Sigma X)$ by $(a_1|\ldots |a_n)$,
where $a_i\in \widetilde \C_*(X)$. Since the diagonal mapping $\Sigma X \to
\Sigma X \times \Sigma X$ gives the $H$-map $\Delta \colon
\Omega\Sigma X \to \Omega\Sigma X \times \Omega\Sigma X$, it follows
from the Eilenberg-Zilber theorem that
\[ \Delta_*(a_1|\ldots |a_n) = (\Delta_* a_1|\ldots |\Delta_* a_n)  \]
where $(a_1\otimes b_1 | a_2\otimes b_2) =
(a_1|a_2)\otimes(b_1|b_2)$.

There is a nice combinatorial model for any topological space of the
form $\Omega\Sigma X$ with $X$ connected, namely the James
construction $JX$ on $X$. After one suspension this gives rise to a
splitting
\[ \Sigma\Omega\Sigma X \sim \Sigma JX \sim\bigvee_{n\geqslant1}\Sigma X^{(n)}, \]
where $X^{(n)}$ denotes the $n$-fold smash power of $X$.

\begin{example} There exists a homotopy equivalence
 \[ \Sigma\Omega\Sigma S^2 \to \bigvee\limits_{n\geqslant1} \Sigma(S^2)^{(n)}\simeq \Sigma
\Big(\bigvee\limits_{n\geqslant1} S^{2n}\Big). \]  Therefore there exists an $H$-map \, $\Omega \Sigma(\Omega\Sigma S^2)
\to \Omega \Sigma \Big(\bigvee\limits_{n\geqslant1} S^{2n}\Big)$ that is
a homotopy equivalence.
\end{example}

Using these classical topological results let us describe topological
realizations of the Hopf algebras we study in this work.

We have the following Hopf algebras:
\begin{itemize}
\item[\textbf{I.}]
\begin{enumerate}
\item $\Coh_*(\mathbb CP^{\infty})$ is a divided power algebra $\mathbb
Z[v_1,v_2,\dots]/I$, where the ideal $I$ is generated by the relations $v_nv_m-{n+m\choose
n}v_{n+m}$, with the comultiplication
\begin{equation}\label{F-1}
\Delta v_n=\sum_{k=0}^{n}v_k\otimes v_{n-k};
\end{equation}
\item $\Coh^*(\mathbb CP^{\infty})=\mathbb Z[u]$ -- a polynomial ring
with the comultiplication
$$
\Delta u=1\otimes u+u\otimes 1;
$$
\end{enumerate}

\item[\textbf{II.}]
\begin{enumerate}
\item $\Coh_*(\Omega\Sigma\mathbb CP^{\infty})\simeq
T(\tilde\Coh_*(\mathbb CP^{\infty}))=\mathbb Z\langle
v_1,v_2,\dots\rangle$ with $v_i$ being non-commuting variables of
degree $2i$. Thus there is an isomorphism of Hopf algebras
$$
\Coh_*(\Omega\Sigma\mathbb CP^{\infty})\simeq\mathcal{Z}
$$
under which $v_n$ corresponds to $Z_n$.

The coproduct $\Delta$ on $H_*(\Omega\Sigma\mathbb CP^{\infty})$
induced by the diagonal in $\Omega\Sigma\mathbb CP^{\infty}$ is
compatible with the one in $\mathcal{Z}$:
$$
\Delta v_n=\sum\limits_{i+j=n}v_i\otimes v_j
$$
So if we set $\deg Z_i=2i$, then there is an isomorphism of graded
Hopf algebras. This gives a geometric interpretation for the
antipode $\chi$  in $\mathcal{Z}$ : in $\Coh_*(\Omega\Sigma\mathbb
CP^{\infty})$ it arises from the time-inversion of loops.

\item $\Coh^*(\Omega\Sigma\mathbb CP^{\infty})$ is the graded dual Hopf
algebra to $\Coh_*(\Omega\Sigma\mathbb CP^{\infty})$.
\end{enumerate}

\begin{thm*}\, [BR]
Let us double the graduation, that is let us set $\deg Z_i=2i$,
$\deg t_j=2j$. Then there is an isomorphism of graded Hopf algebras:
$$
\Coh_*(\Omega\Sigma \mathbb CP^{\infty})\simeq \mathcal{Z}=\mathbb Z \langle
Z_1,Z_2,\dots\rangle,\quad \Coh^*(\Omega\Sigma \mathbb CP^{\infty})\simeq
\mathcal{M}=\Qsym[t_1,t_2,\dots].
$$
\end{thm*}

\item[\textbf{III.}]
$\Coh_*(BU)\simeq \Coh^*(BU)\simeq\mathbb
Z[\sigma_1,\sigma_2,\dots]\simeq\mathcal{C}$. It is a self-dual Hopf
algebra of symmetric functions. In the cohomology $\sigma_i$ are
represented by Chern classes.

\item[\textbf{IV.}]
\begin{enumerate}
\item $\Coh_*(\Omega\Sigma S^2)=\Coh_*(\Omega S^3)=\mathbb Z[v]$
-- a polynomial ring with $\deg v=2$ and the comultiplication
$$
\Delta v=1\otimes v+v\otimes 1;
$$

\item $\Coh^*(\Omega\Sigma S^2)=\mathbb
Z[v_n]/I$ is a divided power algebra. Thus $\Coh^*(\Omega\Sigma
S^2)\simeq\Coh_*(\mathbb CP^{\infty})$.
\end{enumerate}

\item[\textbf{V.}]
$\Coh_*(\Omega\Sigma(\Omega\Sigma S^2))\simeq\mathbb Z\langle
w_1,w_2,\dots\rangle$. It is a free associative Hopf algebra with the
comultiplication
\begin{equation}\label{F-2}
\Delta w_n=\sum\limits_{k=0}^n{n\choose k}w_k\otimes w_{n-k}
\end{equation}

\item[\textbf{VI.}]
$\Coh_*\Big(\Omega\Sigma\Big(\bigvee\limits_{n\geqslant1}S^{2n}\Big)\Big)\simeq
\mathbb{Z}\langle \xi_1,\xi_2,\ldots\rangle$. It is a free
associative algebra and has the structure of a graded Hopf algebra with the comultiplication
\begin{equation}\label{F-3}
\Delta\xi_n = 1\otimes \xi_n + \xi_n\otimes 1.
\end{equation}
Therefore,
$\Coh_*\Big(\Omega\Sigma\Big(\bigvee\limits_{n\geqslant1}S^{2n}\Big)\Big)$
gives a topological realization of the universal Lie-Hopf algebra
$\mathcal{U}$.

The homotopy equivalence
\[ a \colon \Omega\Sigma\Big(\bigvee\limits_{n\geqslant1}S^{2n}\Big)
\longrightarrow \Omega\Sigma(\Omega\Sigma S^2) \] induces
an isomorphism of graded Hopf algebras
 \[ a_* \colon \mathbb{Z}\langle \xi_1,\xi_2,\ldots\rangle \longrightarrow
 \mathbb{Z}\langle w_1,w_2,\ldots\rangle, \]
and its algebraic form is determined by the conditions:
 \[ \Delta a_*\xi_n = (a_* \otimes a_*)(\Delta\xi_n). \]
For example, $a_*\xi_1 = w_1,\; a_*\xi_2 = w_2- w_1|w_1,\; a_*\xi_3 =
w_3- 3w_2|w_1 + 2w_1|w_1|w_1$.
\end{itemize}

Thus using topological results we have obtained that two Hopf algebra structures
on the free associative algebra with the comultiplications  (\ref{F-2}) and (\ref{F-3}) are isomorphic over $\mathbb{Z}$.

This result is interesting from the topological point of view, since the elements $(w_n-a_*\xi_n)$ for $n\geqslant 2$ are obstructions
to the desuspension of the homotopy equivalence

$$
\Sigma(\Omega\Sigma S^2) \to
\Sigma\Big(\bigvee\limits_{n\geqslant1}S^{2n}\Big).
$$

We have the commutative diagram:
$$
\xymatrix{
\Omega\Sigma(\Omega\Sigma S^2)\ar[r]^{\Omega\Sigma k} &
\Omega\Sigma \mathbb{C}P^{\infty}\ar@{-->}[dr]^j&&\\
\Omega\Sigma S^2\ar[r]^{k} \ar[u] & \mathbb{C}P^{\infty}
\ar[u]\ar[r]^i&BU\ar[r]^{\det}&\mathbb{C}P^{\infty}\\
}
$$
Here
\begin{itemize}
\item $i$ is the inclusion $\mathbb{C}P^{\infty}=BU(1)\subset BU$,

\item $j$ arises from the universal property of $\Omega\Sigma
\mathbb{C}P^{\infty}$ as a free $H$-space;

\item the mapping $\mathbb{C}P^{\infty}\to\Omega\Sigma \mathbb{C}P^{\infty}$ is
induced by the identity map $\Sigma \mathbb{C}P^{\infty}\to\Sigma
\mathbb CP^{\infty}$,

\item the mapping $k\colon \Omega\Sigma S^2 \to
\mathbb CP^{\infty}$ corresponds to the generator of
$\Coh^2(\Omega\Sigma S^2)=\mathbb{Z}$.

\item $\det$ is the mapping of the classifying spaces $BU\to
BU(1)$  induced by the mappings \linebreak $\det: U(n)\to
U(1),\;n=1,2,\dots$.
\end{itemize}

As we have mentioned, $\Coh^*(BU)$ as a Hopf algebra is
isomorphic to the algebra of symmetric functions (in the
generators -- Chern classes) and is self-dual.
$$
\Coh^*(BU)\simeq \Coh_*(BU)\simeq\Sym[t_1,t_2,\dots]=\mathbb
Z[\sigma_1,\sigma_2,\dots].
$$
Since $\mathbb CP^{\infty}$ gives rise to the algebra generators in
$\Coh_*(BU)$, $j_*$ is an epimorphism on homology, and $j^*$ is a
monomorphism on cohomology.
\begin{itemize}
\item $j_*$ corresponds to the factorization of the
non-commutative polynomial algebra $\mathcal{Z}$ by the
commutation relations $Z_iZ_j-Z_jZ_i$ and sends $Z_i$ to
$\sigma_i$.

\item $j^*$ corresponds to the inclusion
$\Sym[t_1,t_2,\dots]\subset\Qsym[t_1,t_2,\dots]$ and sends
$\sigma_i$ to $M_{\omega_i}$, where
$\omega_i=\underbrace{(1,1,\dots,1)}_i$.

\item $k^*$ sends $u$ to $v_1$ and defines an embedding of the polynomial ring $\mathbb Z[u]$
into the divided power algebra $\mathbb Z[v_n]/I$.

\item the composition $\det_*\circ j_*$ maps the algebra
$\mathcal{Z}$ to the divided power algebra $H_*(\mathbb CP^{\infty})$.
This corresponds to the mapping $Z_i\to
\left.d_k\right|_{\mathcal{P}_s}$ of $\mathcal{Z}$ to
$\mathcal{D}(\mathcal{P}_s)$.
\end{itemize}

\section{Flag $f$-vectors}
\begin{defin}
Let $P^n$ be an $n$-dimensional polytope and
$S=\{a_1,\dots,a_k\}\subset\{0,1,\dots,n-1\}$.\\
A {\itshape flag
number} $f_{S}=f_{a_1,\,\dots,\,a_k}$ is the number of  increasing sequences of faces
$$
F^{a_1}\subset F^{a_2}\subset\dots\subset F^{a_k},\quad\dim F^{a_i}=a_i.
$$
It is easy to see that for $S=\{i\}$ the number $f_{\{i\}}=f_i$ is just the
number of $i$-dimensional faces. We have already defined this number above.

The collection $\{f_{S}\}$ of all the flag numbers is called a
{\itshape flag $f$-vector} (or {\itshape extended $f$-vector}) of
the polytope $P^n$. By the definition $f_{\varnothing}=1$.
\end{defin}

Flag $f$-vectors have been extensively studied by M.~Bayer and
L.~Billera in \cite{BB}, where the generalized Dehn-Sommerville
relations are proved:
\begin{thm*}
Let $P^n$ be an $n$-dimensional polytope, and
$S\subset\{0,\dots,n-1\}$. If $\{i,k\}\subseteq S\cup\{-1,n\}$
such that $i<k-1$ and $S\cap\{i+1,\dots,k-1\}=\varnothing$, then
$$
\sum\limits_{j=i+1}^{k-1}(-1)^{j-i-1}f_{S\cup\{j\}}=(1-(-1)^{k-i-1})f_S.
$$
\end{thm*}
\begin{defin}
For $n\geqslant 1$ let $\Psi^n$ be the set of subsets
$S\subset\{0,1,\dots,n-2\}$ such that $S$ contains no two
consecutive integers.
\end{defin}
It easy to show by induction that the cardinality of $\Psi^n$ is
equal to the $n$-th Fibonacci number $c_n$ ($c_n=c_{n-1}+c_{n-2}$,
$c_0=1,\;c_1=1$).

For any polytope $P$ there exists a {\itshape cone} (or a
{\itshape pyramid}) $CP$ and a {\itshape suspension} (or a
{\itshape bipyramid}) $BP$. These two operations are defined on
combinatorial polytopes and can be extended to linear operators on
the ring $\mathcal{P}$. It is natural to set
$B\varnothing=1=C\varnothing$.

The face lattice of the polytope $CP^n$ is:
$$
\{\varnothing\};\;\{F^0_i,
C\varnothing\};\;\{F^1_i,CF^0_j\};\;\dots;\;\{F^{n-1}_i,CF^{n-2}_j\},\;\{P^n,CF^{n-1}_i\};\{CP^n\},
$$
where $F_i^k$ are $k$-dimensional faces of $P^n$, and
\begin{gather*}
F^{a_1}\subset F^{a_2}\;\mbox{in $CP$} \;\Leftrightarrow\;
F^{a_1}\subset
F^{a_2}\;\mbox{in $P$},\\
F^{a_1}\subset CF^{a_2}\;\Leftrightarrow \;F^{a_1}\subset
F^{a_2},\\
CF^{a_1}\subset CF^{a_2}\;\Leftrightarrow\;F^{a_1}\subset F^{a_2}.
\end{gather*}

The polytope $BP^n$ has the face lattice:
$$
\{\varnothing\},\{C_0\varnothing, C_1\varnothing
,F^0_i\},\{F^1_i,C_0F^0_j,C_1F^0_k\},\dots,\{F_k^{n-1},C_0F^{n-2}_j,
C_1F^{n-2}_k\},\{C_0F^{n-1}_i,C_1F^{n-1}_j\},\{BP^n\}.
$$
where $C_0$ and $C_1$ are the lower and the upper cones, and
\begin{gather*}
F^{a_1}\subset F^{a_2}\;\mbox{in $BP$}\Leftrightarrow\; F^{a_1}\subset F^{a_2}\;\mbox{in $P$},\\
F^{a_1}\subset C_sF^{a_2}\;\Leftrightarrow F^{a_1}\;\subset F^{a_2};\\
C_sF^{a_1}\subset C_tF^{a_2}\;\Leftrightarrow\; s=t\mbox{
and}\;F^{a_1}\subset F^{a_2};
\end{gather*}
\begin{defin}
For $n\geqslant 1$ let  $\Omega^n$ be the set of $n$-dimensional polytopes that arises when we apply  words in $B$ and $C$ that end in $C^2$ and contain no adjacent $B$'s to the empty set $\varnothing$.
\end{defin}
Each word of length $n+1$ from the set $\Omega^n$ either has the
form $CQ$, $Q\in\Omega^{n-1}$, or $BCQ$, $Q\in\Omega^{n-2}$, so
cardinality of the set $\Omega^n$ satisfies the Fibonacci relation
$|\Omega^n|=|\Omega^{n-1}|+|\Omega^{n-2}|$. Since
$|\Omega^1|=|\{C^2\}|=1$, and $|\Omega^2|=|\{C^3,BC^2\}|=2$, we see
that $|\Omega^n|=c_n=|\Psi^n|$.

M.~Bayer and L.~Billera proved the following fact:
\begin{thm*}
Let $n\geqslant 1$. Then
\begin{enumerate}
\item
For all $T\subseteq\{0,1,\dots,n-1\}$ there is a nontrivial linear
relation expressing $f_T(P)$ in terms of $f_S(P),\;S\in\Psi^n$,
which holds for all $n$-dimensional polytopes (see
\cite[Proposition 2.2]{BB}).
\item The extended $f$-vectors of the $c_n$ elements of $\Omega^n$
are affinely independent \cite[Proposition 2.3]{BB}. Thus the flag
$f$-vectors
$$
\{f_S(P^n)\}_{S\in\Psi^n}:P^n\mbox{ -- an $n$-dimensional
polytope}
$$
span an $(c_n-1)$-dimensional affine hyperplane defined by the
equation $f_{\varnothing}=1$ (see \cite[Theorem 2.6]{BB})
\end{enumerate}
\end{thm*}

In fact, a little bit more stronger fact which can be easily
extracted from the original Bayer-Billera's proof is true.

Let us identify the words in $\Omega^n$ with the sets in $\Psi^n$
in such a way that the word \linebreak
$C^{n+1-a_k}BC^{a_k-a_{k-1}-1}B\dots BC^{a_1-1}$ corresponds to
the set $\{a_1-3,\dots,a_k-3\}$. Let us set $C<B$ and order the
words lexicographically. Consider a matrix $K^n$ of sizes
$c_n\times c_n$, $k_{Q,\,S}=f_{S}(Q),\;Q\in\Omega^n,\,S\in \Psi^n,
$.
\begin{prop}
$$
\det(K^n)=1.
$$
\end{prop}
\begin{proof}
For $n=1$ the matrix $K^1=(f_{\varnothing}(I))=(1)$.

For $n=2$
$$
K^2=\begin{pmatrix}
1&f_{0}(\Delta^2)\\
1&f_{0}(I^2)\\
\end{pmatrix}=
\begin{pmatrix}
1&3\\
1&4\\
\end{pmatrix};\quad\det K^2=1
$$
Let us prove the statement by induction. Consider the matrix
$K^n$:
$$
K^n=\begin{pmatrix}
K_{11}& K_{12}\\
K_{21}& K_{22}\\
\end{pmatrix}
$$
where the block $K_{11}$ corresponds to the words of the form
$CQ,\;Q\in\Omega^{n-1}$ and the sets $S$, that do not contain
$n-2$. Similarly, the block $K_{22}$ corresponds to the words of
the form $BCQ,\;Q\in\Omega^{n-2}$ and the sets containing $n-2$.

Each increasing sequence of faces of the polytope $CP^{n-1}$
without two faces of adjacent dimensions has the form
$$
F^{l_1}\subset\dots\subset F^{l_i}\subset
CF^{l_{i+1}}\subset\dots\subset CF^{l_k},
$$
so for $\{a_1,\dots,a_k\}\in\Psi^n$
\begin{equation}\label{X2}
f_{a_1,\,\dots,\,a_k}(CP^{n-1})=f_{a_1-1,\,a_2-1,\,\dots,\,a_k-1}(P^n)+\dots+f_{a_1,\,\dots,\,a_i,\,a_{i+1}-1,\,\dots,\,a_k-1}(P^n)+\dots+f_{a_1,\,\dots,\,a_k}(P^n).
\end{equation}
where $f_{-1,\,a_2,\,\dots,\,a_k}=f_{a_2,\,\dots,\,a_k}$.

The sets of the form $\{a_1,\dots,a_l,a_{l+1}-1,\dots,a_k-1\}$ may
not belong to $\Psi^{n-1}$. But we can express them in terms of the
sets from $\Psi^{n-1}$ using the generalized Dehn-Sommerville
relations. Let $a_k<n-2$. Each relation has the form:
$$
f_{S_1,\,i,\,k-1,\,k,\,S_2}=(-1)^{k-i}\left(\left(1-(-1)^{k-i-1}\right)f_{S_1,\,i,\,k,\,S_2}-\sum\limits_{j=i+1}^{k-2}(-1)^{j-i-1}f_{S_1,\,i,\,j,\,k,\,S_2}\right),
$$
where all the sets on the right side are lower than the set on the
left side. Thus for $a_k<n-2$
$$
f_{a_1\,\dots,\,a_k}(CP^n)=f_{a_1,\,\dots,\,a_k}(P^n)+\mbox{lower
summands},
$$
so the matrix $K_{11}$ can be represented as $K^{n-1}T$, where $T$
is an upper unitriangular matrix. In particular, $\det K_{11}=\det
K^{n-1}\cdot\det T=1$ by induction.

\begin{lemma}
Let $P$ be an $(n-2)$-dimensional polytope, and $S\in\Psi^n$. Then
for the polytopes $BCP$ and $CBP$
$$
f_{S}(BCP)=
\begin{cases}
f_{S}(CBP),&\mbox{if $n-2\notin S$},\\
f_{S}(CBP)+f_{S\setminus\{n-2\}}(P),&\mbox{if $n-2\in S$},
\end{cases}\label{BC-CB}
$$
\end{lemma}
\begin{proof}
For the polytope $BCP$ each increasing sequence of faces has one
of the forms
\begin{gather*}
F^{l_1}\subset\dots \subset F^{l_i}\subset CF^{l_{i+1}}\subset
\dots\subset CF^{l_j}\subset C_sCF^{l_{j+1}}\subset\dots\subset
C_sCF^{l_k}\;\mbox{or}\\
F^{l_1}\subset\dots\subset F^{l_i}\subset C_sF^{l_{i+1}}\subset
\dots\subset C_sF^{l_j}\subset C_sCF^{l_{j+1}}\subset\dots\subset
C_sCF^{l_k}\;\mbox{or}\\
F^{l_1}\subset\dots\subset F^{l_j}\subset
C_sCF^{l_{j+1}}\subset\dots\subset C_sCF^{l_k},\;\mbox{when
$i=j$},
\end{gather*}
where $s=0$ or $1$.

Similarly, for the polytope $CBP$ each increasing sequence of
faces has one of the forms
\begin{gather*}
F^{l_1}\subset\dots \subset F^{l_i}\subset CF^{l_{i+1}}\subset
\dots\subset CF^{l_j}\subset CC_sF^{l_{j+1}}\subset\dots\subset
CC_sF^{l_k}\;\mbox{or}\\
F^{l_1}\subset\dots\subset F^{l_i}\subset C_sF^{l_{i+1}}\subset
\dots\subset C_sF^{l_j}\subset CC_sF^{l_{j+1}}\subset\dots\subset
CC_sF^{l_k}\;\mbox{or}\\
F^{l_1}\subset\dots\subset F^{l_j}\subset
CC_sF^{l_{j+1}}\subset\dots\subset CC_sF^{l_k},\;\mbox{when
$i=j$},
\end{gather*}
where $s=0$ or $1$, except for the additional condition
$F^{l_k}\ne P$, that is $l_k\ne n-2$.

Thus for $l_k\ne n-2$, we can exchange $CC_s$ and $C_sC$ to obtain
a one-to-one correspondence between the increasing sequences. So
$f_S(BCP)=f_S(CBP)$, if $n-2\notin S$.

If $l_k=n-2$, then we can exchange $CC_s$ and $C_sC$ as before except for the case when the sequence of faces in the polytope $BCP$
has the form
$$
F^{a_1}\subset\dots\subset F^{a_{k-1}}\subset P
$$
The sequences of this type give exactly
$f_{a_1,\,\dots,\,a_{k-1}}(P^{n-1})$, so
$f_{S}(BCP)=f_{S}(CBP)+f_{S\setminus\{n-2\}}(P)$, if $n-2\in S$.
\end{proof}

Let us consider the polytope $BCQ$ corresponding to one of the
lower rows of $K^n$. If $Q$ starts with $C$, then there is a row
$CBQ$ in the upper part of the matrix. Let us subtract the row
$CBQ$ from the row $BCQ$. Then lemma \ref{BC-CB} implies that the
resulting row is
\begin{equation}\label{X3}
k'_{BCQ,\,S}=
\begin{cases}
0,& \mbox{if }n-2\in S;\\
f_S(Q),&\mbox{if } n-2\in S.
\end{cases}
\end{equation}

If $Q$ starts with $B$, then there is no row $CBQ$ in the matrix.
But since $\det K^{n-1}=1$ and all flag $f$-numbers of
$(n-1)$-dimensional polytopes can be expressed in terms of
$\{f_{S},\;S\in\Psi^{n-1}\}$, flag $f$-vector of any
$(n-1)$-dimensional polytope is an integer combination of flag
$f$-vectors of the polytopes from $\Omega^{n-1}$. So flag
$f$-vector of the polytope $BQ$ can be expressed as an integer
combination of flag $f$-vectors of the polytopes
$Q'\in\Omega^{n-1}$:
$$
f(BQ)=\sum\limits_{Q'}n_{Q'}f(Q')
$$
Using the formula (\ref{X2}) we obtain:
$$
f(CBQ)=\sum\limits_{Q'}n_{Q'}f(CQ')
$$
If we subtract the corresponding integer combination of the rows
of the upper part of the matrix from the row $BCQ$, we obtain the
row (\ref{X3}).

Thus we see, that using elementary transformations of rows the
matrix $K^{n-1}$ can be transformed to the matrix
$$
\begin{pmatrix}
K^{n-1}T&K_{12}\\
0&K^{n-2}
\end{pmatrix}
$$
By the inductive assumption $\det K^{n-1}=\det K^{n-2}=1$, so
$\det K^n=1$.
\end{proof}
\begin{rem}
In the proof we follow the original Bayer-Billera's idea, except
for the fact that they use the additional matrix of face numbers
of the polytopes $Q\in\Psi^n$, and our proof is direct.
\end{rem}
\begin{cor}
Flag $f$-vector of any $n$-dimensional polytope $P^n$ is an
integer combination of flag $f$-vectors of the polytopes
$Q\in\Omega^n$.
\end{cor}
\begin{proof}
Indeed, any flag $f$-number is a linear combination of
$f_S,\;S\in\Psi^n$. Since $\det K^n=1$, the vector
$\{f_S(P),\;S\in\Psi^n\}$ is an integer combination of the vectors
$\{f_S(Q),\;S\in\Psi^n\},\;Q\in\Omega^n$.

This implies that the whole vector $\{f_S(P)\}$ is an integer
combination of the vectors $\{f_S(Q)\},\;Q\in\Omega^n$ with the
same coefficients.
\end{proof}
This fact is  important, when we try to describe an image of the
generalized $f$-polynomial (see below) in the ring
$\Qsym[t_1,t_2,\dots][\alpha]$.

Generalized $f$-vectors are also connected with a very interesting
construction of the $cd$-index invented by J.Fine (see \cite[Prop.
2]{BK}, see also the papers by R.~Stanley \cite{St3} and M.~Bayer,
A.~Klapper \cite{BK}).

\section{Operators $B$ and $C$}
As it was mentioned below, the operations $C$ and $B$ can be
extended to the linear operators on the ring of polytopes.

Then for any $n$-dimensional polytope $P^n$ we have:
\begin{gather*}
d_kCP^n=Cd_kP^n+d_{k-1}P^n,\; 1\leqslant k\leqslant n,\quad
d_{n+1}CP^n=1+d_nP^n=C\varnothing+d_nP^n;\\
dBP^n=2CP^n,\;n\geqslant 1;\quad dB1=dI=2=2C\varnothing,\\
d_kBP^n=2Cd_kP^n+d_{k-1}P^n,\;2\leqslant k\leqslant n,\quad
d_{n+1}BP^n=2+d_nP^n=2C\varnothing+d_nP^n.
\end{gather*}
Thus if we add to the ring $\mathcal{P}$ the new element
$\varnothing$ of the graduation $-1$ such that
$\varnothing=d1,\;d\varnothing=0$, then the operators $B$ and $C$
satisfy the relations:
\begin{gather*}
dC-Cd=[d,C]=1;\quad[d_k,C]=d_{k-1},\;k\geqslant 2;\\
[d,B]=(2C-B)d;\quad[d_k,B]=(2C-B)d_k+d_{k-1};\;k\geqslant 2.
\end{gather*}
Otherwise the relations in the ring $\mathcal{P}$ have the form:
\begin{gather*}
[d,C]=1+\xi_0;\quad[d_{k+1},C]=d_k+\frac{1}{k!}\left.\frac{\partial^k}{\partial\alpha^k}\right|_{\alpha=0}\xi_{\alpha},\;k\geqslant
1;\\
[d,B]=(2C-B)d+2\xi_0;\quad[d_{k+1},B]=(2C-B)d_k+d_{k-1}+\frac{2}{k!}\left.\frac{\partial^k}{\partial\alpha^k}\right|_{\alpha=0}\xi_{\alpha},\;k\geqslant
1;
\end{gather*}
It is convenient to write these relations in terms of the operator
$\Phi(t)$. Let us denote $A=2C-B$. Then
\begin{gather*}
[\Phi(t),C]=t\Phi(t)+t\xi_t;\\
[\Phi(t),B]=(A+t)\left(\Phi(t)-1\right)+2t\xi_t;\\
\Phi(t)A=A+t+t\Phi(t);\\
\Phi(t)[B,C]=[B,C]+tA+t^2;
\end{gather*}
\begin{prop}
For any $m\geqslant 0$ the formulas:
\begin{gather*}
(Ag)(\alpha,t_1,\dots,t_{m+1})=(\alpha+t_{m+1})g(\alpha,0,\dots,0)+(t_m+t_{m+1})g(\alpha,t_{m+1},0,\dots,0,)+\dots+\\
+(t_{i-1}+t_i)g(\alpha,t_i,t_{i+1},\dots,t_{m+1},0,\dots,0,)+\dots+t_1(v_{m+1}g)(\alpha,t_1,\dots,t_{m+1});\\
(Cg)(\alpha,t_1,\dots,t_{m+1})=(\alpha+t_1+\dots+t_{m+1})(v_{m+1}g)(\alpha,t_1,\dots,t_{m+1})+t_{m+1}g(t_{m+1},t_1,\dots,t_m)+\\+
t_mg(t_m,t_1,\dots,t_{m-1},0)+\dots+t_ig(t_i,t_1,\dots,t_{i-1},0,\dots,0)+\dots+t_1g(t_1,0,\dots,0)
\end{gather*}
define the linear mappings
$A,C:\Qsym[t_1,\dots,t_m][\alpha]\to\Qsym[t_1,\dots,t_{m+1}][\alpha]$
such that
$$
f_{m+1}(AP^n)=(Af_m)(P^n), \mbox{ and } f_{m+1}(CP^n)=(Cf_m)(P^n)
$$
for any $P^n\in\mathcal{P}^n$ and $m\geqslant n$.
\end{prop}
\begin{proof}
\begin{multline*}
\Phi(t_{m+1})\Phi(t_m)\dots\Phi(t_1)AP=\Phi(t_{m+1})\Phi(t_m)\dots\Phi(t_2)\left(A+t_1\Phi(t_1)+t_1\right)P=\\
=\Phi(t_{m+1})\Phi(t_m)\dots\Phi(t_3)\left(A+t_2+(t_2+t_1)\Phi(t_2)+t_1\Phi(t_2)\Phi(t_1)\right)P=\\
=\left(A+t_{m+1}+(t_{m+1}+t_m)\Phi(t_{m+1})+(t_m+t_{m-1})\Phi(t_{m+1})\Phi(t_m)+\dots+t_1\Phi(t_{m+1})\dots\Phi(t_1)\right)P
\end{multline*}
This gives the formula for $f_{m+1}(AP^n)$, and for $m\geqslant n$
$$
\xi_{\alpha}\Phi(t_{m+1})\dots\Phi(t_1)P^n=f_{m+1}(\alpha,t_1,\dots,t_{m+1})(P^n)=(v_{m+1}f_m)(\alpha,t_1,\dots,t_{m+1})(P^n)
$$
Similarly we have:
\begin{multline*}
\Phi(t_{m+1})\Phi(t_m)\dots\Phi(t_1)CP=\Phi(t_{m+1})\dots\Phi(t_2)\left((C+t_1)\Phi(t_1)+t_1\xi_{t_1}\right)P=\\
=\Phi(t_{m+1})\Phi(t_m)\dots\Phi(t_3)\left((C+t_2+t_1)\Phi(t_2)\Phi(t_2)+t_2\xi_{t_2}\Phi(t_1)+t_1\xi_{t_1}\right)P=\\
\left((C+t_{m+1}+\dots+t_1)\Phi(t_{m+1})\dots\Phi(t_1)P+t_{m+1}\xi_{t_{m+1}}\Phi(t_m)\dots\Phi(t_1)+\dots+t_1\xi_{t_1}\right)P=\\
=(C+t_{m+1}+\dots+t_1)\Phi(t_{m+1})\dots\Phi(t_1)P+t_{m+1}f_m(t_{m+1},t_1,\dots,t_m)(P)+\\
+t_mf_m(t_m,t_1,\dots,t_{m-1},0)(P)+\dots+t_1f(t_1,0,\dots,0)(P)
\end{multline*}

Let us prove that for a polynomial $g\in\Qsym[t_1,\dots,t_m][\alpha]$ the polynomials $Ag$ and $Cg$ belong to $\Qsym[t_1,\dots,t_{m+1}][\alpha]$.
\begin{multline*}
(Ag)(\alpha,t_1,\dots,t_{i-1},0,t_{i+1},\dots,t_{m+1})=(\alpha+t_{m+1})g(\alpha,0,\dots,0)+\dots+(t_{i+1}+0)g(\alpha,t_{i+1},\dots,t_{m+1},0,\dots,0)+\\
+(0+t_{i-1})g(\alpha,0,t_{i+1},\dots,t_{m+1},0,\dots,0)+\dots+t_1g(\alpha,t_1,\dots,t_{i-1},t_{i+1},t_m)=\\
=(\alpha+t_{m+1})g(\alpha,0,\dots,0)+\dots+(t_{i+1}+t_{i-1})g(\alpha,t_{i+1},\dots,t_{m+1},0,\dots,0)+\dots+t_{m+1}g(\alpha,t_1,\dots,t_{i-1},t_{i+1},\dots,0)=\\
=(Ag)(\alpha,t_1,\dots,t_{i-1},t_{i+1},\dots,t_{m+1},0)
\end{multline*}
The same situation we have for $C$:
\begin{multline*}
(Cg)(\alpha,t_1,\dots,t_{i-1},0,t_{i+1},\dots,t_{m+1})=(\alpha+\dots+t_{i-1}+t_{i+1}+\dots+t_{m+1})g(t_{m+1},t_1,\dots,t_{i-1},t_{i+1},\dots,t_{m+1})+\\
+\dots+t_{i+1}g(t_{i+1},t_1,\dots,t_{i-1},0,\dots,0)+0g(0,t_1,\dots,t_{i-1},0,\dots,0)+t_{i-1}g(t_{i-1},t_1,\dots,t_{i-2},0,\dots,0)+\\
+\dots+t_1g(t_1,0,\dots,0)=(Cg)(\alpha,t_1,\dots,t_{i-1},t_{i+1},\dots,t_{m+1},0)
\end{multline*}
\end{proof}
The operator $B$ can be defined as $B=2C-A$. Then
$f_{m+1}(BP^n)=(Bf_m)(P^n)$ for $m\geqslant n$.

For the same reasons the operators $A, B, C$ are defined on $\Qsym[t_1,t_2,\dots][\alpha]$:
\begin{gather*}
(Ag)(\alpha,t_1,t_2,\dots)=\alpha g(\alpha,0,\dots,)+t_1g(\alpha,t_1,t_2,\dots)+\sum\limits_{i=2}^{\infty}(t_{i-1}+t_i)g(\alpha,t_i,t_{i+1},\dots);\\
(Cg)(\alpha,t_1,t_2,\dots)=(\alpha+\sigma_1)g(\alpha,t_1,t_2,\dots)+\sum\limits_{i=1}^{\infty}t_ig(t_i,t_1,\dots,t_{i-1},0,\dots,0);\\
B=2C-A,
\end{gather*}
and $(Af)(P)=f(AP),\;(Bf)(P)=f(BP),\;(Cf)(P)=f(CP)$.

It is easy to see that
\begin{gather*}
A\left(\alpha^mM_{(a_1,\,\dots,\,a_k)}\right)=\alpha^m\left(\alpha+2M_{(1,\,a_1,\,\dots,\,a_k)}+M_{(a_1+1,\,a_2,\,\dots,\,a_k)}\right),\\
C\left(\alpha^mM_{(a_1,\,\dots,\,a_k)}\right)=(\alpha+\sigma_1)M_{(a_1,\,\dots,\,a_k)}+M_{(a_1,\,\dots,\,a_k,\,m+1)}.
\end{gather*}
\section{Structure of $\mathcal{D}$}
\begin{prop}
$$
\Phi(-t)\Phi(t)=1,
$$
that is for any $n\geqslant 1$
\begin{equation}\label{X4}
d_n-d_1d_{n-1}+\dots+(-1)^{n-1}d_{n-1}d_1+(-1)^nd_n=0
\end{equation}
\end{prop}
To prove this fact we need the notion of a {\itshape face figure}.
\begin{defin}
Let $F$ be an $i$-face of the $n$-dimensional polytope $P$. A
{\itshape face figure} $P/F$ is defined as
$(F^{\diamondsuit})^*$.
\end{defin}
The face figure $P/F$ is an $(n-i-1)$-dimensional polytope.
$(j-i-1)$-dimensional faces of $P/F$ are in one-to-one
correspondence with $j$-faces $G$ of $P$ such that $F\subset
G\subset P$. Let us denote by $G/F$ the face of $P/F$
corresponding to the face $G$ of $P$.
\begin{proof}
$$
\Phi(-t)\Phi(t)=\sum\limits_{k=0}^{\infty}\left(\sum\limits_{i=0}^k(-1)^id_id_{k-i}\right)t^k.
$$
Let $P^n$ be an $n$-dimensional polytope. If $k>n$ then the
coefficient of $t^k$ in the series $\Phi(-t)\Phi(t)P^n$ is equal
to $0$. If $k=0$, then it is equal to $P^n$. Let $1\leqslant
k\leqslant n$. Then
\begin{multline*}
\left(dd_{k-1}-d_2d_{k-2}+\dots+(-1)^{k-2}d_{k-1}d\right)P^n=\\
=\sum\limits_{F^{n-k}\subset
P^n}\left(\sum\limits_{F^{n-k+1}\supset
F^{n-k}}1-\sum\limits_{F^{n-k+2}\supset
F_{n-k}}1+\dots+(-1)^{k-2}\sum\limits_{F^{n-1}\supset
F^{n-k}}1\right)F^{n-k}=\\
=\sum\limits_{F^{n-k}\subset
P^n}\left(f_0(P^n/F^{n-k})-f_1(P^n/F^{n-k})+\dots+(-1)^{k-2}f_{k-2}(P^n/F^{n-k})\right)F^{n-k}
\end{multline*}
Since  $P^n/F^{n-k}$ is a $(k-1)$-dimensional polytope, the Euler
formula gives the relation
$$
f_0(P^n/F^{n-k})-f_1(P^n/F^{n-k})+\dots+(-1)^{k-2}f_{k-2}(P^n/F^{n-k})=(1+(-1)^k).
$$
Thus we obtain
$$
\left(dd_{k-1}-d_2d_{k-2}+\dots+(-1)^{k-2}d_{k-1}d\right)P^n=(1+(-1)^k)d_kP^n\quad\Leftrightarrow\quad\sum\limits_{i=0}^k(-1)^id_id_{k-i}P^n=0
$$
\end{proof}
\begin{thm}
The homomorphism $R:\mathcal{Z}\to\mathcal{D}$ induces an isomorphism of Hopf algebras:
$$
\mathcal{D}\simeq\mathcal{Z}/J_{\mathcal{D}},
$$
where $J_{\mathcal{D}}$ is a two-sided Hopf ideal
in the Hopf algebra $\mathcal{Z}$, generated by the relations
corresponding to the equality $\Phi(-t)\Phi(t)=1$. \label{D}
\end{thm}
\begin{proof}
Let us prove the following lemma:
\begin{lemma}
Let $D\in\mathcal{D}$ be an operator of graduation $k$. Then for
each space $\mathcal{P}^{[n]},\;n\geqslant k$ there is a unique
representation
\begin{equation}\label{X5}
D=u(d_2,d_3,\dots,d_k)+dv(d_2,d_3,\dots,d_{k-1}),
\end{equation}
and there is a unique representation
\begin{equation}\label{X6}
D=u'(d_2,d_3,\dots,d_k)+v'(d_2,\dots,d_{k-1})d,
\end{equation}
where $u,u',v,v'$ are polynomials.\label{basis}
\end{lemma}
\begin{proof}
For $k=0$ and $1$ it is evident.

Since $2d_2=d^2$, this is true for $k=2$. Let $k\geqslant 3$.

Using the relations (\ref{X4}) we obtain
$$
dd_i=(-1)^id_id+\left(d_2d_{i-1}-d_3d_{i-2}+\dots+(-1)^{i-1}d_{i-1}d_2\right)+(1+(-1)^{i+1})d_{i+1}
$$
So the expressions (\ref{X5}) and (\ref{X6}) exist.

Let $u(d_2,d_3,\dots,d_k)+dv(d_2,\dots,d_{k-1})=0$, where
$u=\sum\limits_{|\omega|=k}a_{\omega}D_{\omega}$, and
$v=\sum\limits_{|\omega|=k-1}b_{\omega}D_{\omega}$.

Then for any $n$-dimensional polytope $P^n,\;n\geqslant k$ we have
$$
\xi_1(u+dv)P^n=0.
$$
Since $u+dv$ has degree $k$, this equality can be written as
$$
\sum\limits_{|\omega|=k}a_{\omega}f_{n-k,\,n-k+j_1,\,n-k+j_1+j_2,\,\dots,\,n-j_l}(P^n)+\sum\limits_{|\omega|=k-1}b_{\omega}f_{n-k,\,n-k+1,\,n-k+1+j_1,\,\dots,\,n-j_l}(P^n).
$$
Using the generalized Dehn-Sommerville equations we obtain
\begin{multline*}
f_{n-k,\,n-k+1,\,n-k+1+j_1,\,\dots,\,n-j_l}=(-1)^{n-k-1}\left(\sum\limits_{j=0}^{n-k-1}(-1)^jf_{j,\,n-k+1,\,n-k+1+j_1,\,\dots,\,n-j_l}\right)+\\
+(1+(-1)^{n-k})f_{n-k+1,\,n-k+1+j_1,\,\dots,\,n-j_l}
\end{multline*}
All the sets $\{j,n-k+1,n-k+j_1,\dots,n-j_l\}$ and
$\{n-k+1,n-k+1+j_1,\dots,n-j_l\}$ on the right are different for
different $\omega$ and belong to $\Psi^n$.

Since the vectors $\{f_{S}(Q),\,S\in\Psi^n\}, Q\in\Omega^n$ are
linearly independent, we obtain that all $a_{\omega}$ and
$b_{\omega}$ are equal to $0$, so the representation (\ref{X5}) is
unique.

We obtain that the monomials $D_{\omega}=d_{j_1}\dots
d_{j_l},\;|\omega|=k,\;j_i\geqslant 2$ and
$dD_{\omega}=dd_{j_1}\dots d_{j_l}\;|\omega|=k-1,\;j_i\geqslant2$
form a basis of the abelian group $\mathcal{D}^k(\mathcal{P}^{[n]})$.

But each monomial $dD_{\omega},\;|\omega|=k,\;j_i\geqslant 2$ can
be expressed as an integer combination of the monomials
$D_{\omega'}$ and $D_{\omega'}d$. So the monomials $D_{\omega}$
and $D_{\omega}d$ also form a basis. This proves the second part
of the lemma.
\end{proof}
Now let us prove the theorem. The mapping
$R:\mathcal{Z}\to\mathcal{D}$ is an epimorphism of Hopf algebras.
Let $z\in\mathcal{Z}$ such that $\deg z=k$ and $Rz=0$.
$$
z=\sum\limits_{|\omega|=k}a_{\omega}Z^{\omega}.
$$
We know that
$$
Z_{k-1}Z_1=(-1)^{k-1}Z_1Z_{k-1}+\sum\limits_{i=2}^{k-2}(-1)^iZ_{k-i}Z_i+(1+(-1)^k)Z_k-\sum\limits_{i=0}^{k}(-1)^iZ_{k-i}Z_i.
$$
So
$$
z=\sum\limits_{|\omega|=k,\;j_i\geqslant 2}a'_{\omega}Z^{\omega}+\sum\limits_{|\omega|=k-1,\;j_i\geqslant 2}b_{\omega}Z_1Z^{\omega}+z',
$$
where $z'\in J_{\mathcal{D}}$.
Since $Rz=Rz'=0$, we obtain
$$
\sum\limits_{|\omega|=k,\;j_i\geqslant 2}a'_{\omega}D_{\omega^*}+\sum\limits_{|\omega|=k-1,\;j_i\geqslant 2}b_{\omega}D_{\omega^*}d=0,
$$
Lemma \ref{basis} implies that all the coefficients $a_{\omega}$ and $b_{\omega}$ are equal to $0$, so $z=z'\in J_{\mathcal{D}}$.
This proves that $\mathcal{D}\simeq \mathcal{Z}/J_{\mathcal{D}}$.
\end{proof}

It is evident that $\left.D_{\omega}\right|_{\mathcal{P}^{[n]}}=0$ for $|\omega|>n$.  The following corollary says that the other relations between operators in $\mathcal{D}$ on the abelian group $\mathcal{P}^{[n]}$ are the same as on the whole ring $\mathcal{P}$.
\begin{cor}
We have
$$
\mathcal{D}(\mathcal{P}^{[n]})=\mathcal{D}/J_n
$$
where the ideal $J_n$ is generated by all the operators $D_{\omega},\;|\omega|>n$.
\end{cor}
\begin{proof}
As in the proof of Theorem \ref{D} using Lemma \ref{basis} it is easy to show that $\mathcal{D}(\mathcal{P}^{[n]})\simeq\mathcal{Z}/J_{\mathcal{D}(\mathcal{P}^{[n]})}$, where the ideal
$J_{\mathcal{D}(\mathcal{P}^{[n]})}$ is generated by the the relations $\Phi(t)\Phi(-t)=1$ and $Z_{\omega}=0$ for $|\omega|>n$.
\end{proof}
\begin{cor}
The operators $d_2,d_3,d_4,\dots$ are algebraically independent.
\end{cor}

\begin{cor}
Rank of the $n$-th graded component of the ring $\mathcal{D}$ is equal
to the $(n-1)$-th Fibonacci number $c_{n-1}$.\label{Fib}
\end{cor}
\begin{proof}
At first, let us calculate rank $r_n$ of the $n$-th graded component of the ring $\mathbb Z\langle d_2,d_3,d_4,\dots\rangle$.
$$
r_0=1,\;r_1=0,\;r_2=1,\;r_3=1,\;r_4=2,\;r_5=3,\dots
$$
It is easy to see that there is a recursive relation
$$
r_{n+1}=r_{n-1}+r_{n-2}+\dots+r_2+1,\;n\geqslant 3
$$
Then $r_{n+1}=r_{n-1}+r_n$. Since $r_2=r_3=1$ we obtain $r_n=c_{n-2},\; n\geqslant 2$.
At last, rank of the $n$-th graded component of the ring $\mathcal{D}$ is equal to $r_{n-1}+r_n=c_{n-3}+c_{n-2}=c_{n-1}$ for $n\geqslant 3$. It is easy to see that for $n=1$ and $2$ it is also true.
\end{proof}

We claim that the antipode in $\mathcal{D}$ is $\chi:d_k\to (-1)^kd_k$. Indeed, the relations
$\Phi(-t)\Phi(t)=1$ exactly mean that $\chi$ is an antipode.
Therefore we obtain the corollary:
\begin{cor}
The Hopf algebra $\mathcal{D}$ is a universal Hopf algebra in the
category of Leibnitz-Hopf algebras with an antipode satisfying the
formula $\chi(H_i)=(-1)^iH_i$.
\end{cor}
\begin{proof}
In fact, let us consider the equation (\ref{X1}) with
$\chi(Z_i)=(-1)^iZ_i$:
$$
(-1)^{n+1}Z_{n+1}+(-1)^nZ_nZ_1+\dots+(-1)Z_1Z_n+Z_{n+1}=0.
$$
The collection of these equalities is equivalent to the equation
$$
\Phi(-t)\Phi(t)=1;
$$
\end{proof}
Let us remind that on the ring of simple polytopes
$\left.d_k\right|_{\mathcal{P}_s}=\left.\frac{d^k}{k!}\right|_{\mathcal{P}_s}$,
so the ring $\mathcal{D}(\mathcal{P}_s)$ is isomorphic to the divided power
ring $d_kd_l={k+l\choose k}d_{k+l}$.
\begin{prop}
$$
\mathcal{D}\otimes\mathbb Q=\mathbb Q\langle
d_1,d_3,d_5,\dots\rangle,
$$
The inclusion $\mathcal{D}\subset\mathcal{D}\otimes\mathbb Q$ is
an embedding, and the operators $d_{2k}$ are expressed in terms of
the operators $d_1,d_3,\dots,d_{2k-1}$ by the formulas
\begin{equation}\label{X7}
d_{2k}=\sum\limits_{i=1}^{k}(-1)^{i-1}\frac{{2i-2\choose i-1}}{i2^{2i-1}}\left(\sum\limits_{j_1+j_2+\dots+j_{2i}=i+k}d_{2j_1-1}\dots
d_{2j_{2i}-1}\right).
\end{equation}
\end{prop}
\begin{proof}
Let us denote
$$
a(t)=\sum\limits_{k\geqslant 0}d_{2k}t^{2k},\quad b(t)=\sum\limits_{k\geqslant 0}d_{2k+1}t^{2k+1}.
$$
Then
\begin{gather*}
\Phi(t)=a(t)+b(t)\;\quad\Phi(-t)=a(t)-b(t);\\
a(t)b(t)=\frac{\Phi(t)+\Phi(-t)}{2}\frac{\Phi(t)-\Phi(-t)}{2}=\frac{\Phi(t)^2-\Phi(-t)^2}{4}=b(t)a(t).
\end{gather*}
Then the relation $\Phi(-t)\Phi(t)=1$ is equivalent to the relations
\begin{gather*}
a(t)^2-b(t)^2=1;\\
a(t)b(t)=b(t)a(t).
\end{gather*}
Therefore $a(t)=\sqrt{1+b(t)^2}$ and the formula (\ref{X7}) is true.
Consequently all the operators $d_{2k}$ are expressed as polynomials in $d_1,d_3,\dots$ with rational coefficients.
For example,
$$
d_2=\frac{d^2}{2},\quad d_4=\frac{dd_3+d_3d}{2}-\frac{d^4}{8}.
$$
This means that the algebra $\mathcal{D}\otimes\mathbb Q$ is generated by $d_1,d_3,d_5,\dots$. On the other hand, let us calculate the number of the monomials
$d_{2j_1-1}\dots d_{2j_k-1}$ with $(2j_1-1)+\dots+(2j_k-1)=n$. Denote this number by $l_n$. Then
$$
l_0=1,\;l_1=1,\;l_2=1,\;l_3=2,\;l_4=3,\dots,
$$
and there is a recursive formula $l_{n+1}=l_{n}+l_{n-2}+l_{n-4}+\dots$ for $n\geqslant 2$. Hence $l_{n+1}=l_n+l_{n-1}$ for $n\geqslant 2$, and $l_n=c_{n-1}$ for $n\geqslant 1$.

We see that the number of monomials of degree $n$ is equal to dimension of the $n$-th graded component according to Corollary \ref{Fib}.
This implies that they are linearly independent over the rationals. Therefore $d_1,d_3,d_5,\dots$ are algebraically independent.
\end{proof}
\begin{defin}
Let us define operators $s_k$ by the formula
$$
s(t)=s_1t+s_2t^2+s_3t^3+\dots=\log\Phi(t).
$$
\end{defin}
Then $s_1=d_1,\;s_2=0,\;s_3=d_3-\frac{d^3}{6}$, and so on.

The relation $\Phi(-t)\Phi(t)=1$ turns into the relation
$s(-t)+s(t)=0$. Thus $s_{2k}=0$ for all $k$.

Also we have $\Delta s(t)=1\otimes s(t)+s(t)\otimes 1$, so each
operator $s_{2k-1}$ is a derivation.
\begin{prop}
There is an isomorphism of Hopf algebras
$$
\mathcal{D}\otimes\mathbb Q=\mathbb Q\langle
s_1,s_3,s_5,\dots\rangle,
$$
where $\mathbb Q\langle s_1,s_3,s_5,\dots\rangle$ is a free Lie-Hopf
algebra in the generators of odd degree, the comultiplication
$\Delta s_{2k-1}=1\otimes s_{2k-1}+s_{2k-1}\otimes 1$, and the antipode $\chi(s_{2k-1})=-s_{2k-1}$\label{S}
\end{prop}

\section{Generalized $f$-polynomial}
The mapping $\Phi:\mathcal{P}\to\mathcal{P}[t]$ is a ring
homomorphism. If we set $\Phi(t_2)t_1=t_1$, then the operator
$\Phi(t_2)$ defines a ring homomorphism
$\mathcal{P}\to\mathcal{P}[t_1,t_2]$: $P\to\Phi(t_2)\Phi(t_1)P$,

Thus for each $n\geqslant 0$ we can define the ring homomorphism
$\Phi_n:\mathcal{P}\to\mathcal{P}[[t_1,\dots,t_n]]$ as a composition:
$$
\Phi_n(t_1,\dots,t_n)P=\Phi(t_n)\dots\Phi(t_1)P.
$$
\begin{defin}
Let us define the mapping
$\varphi:\mathcal{P}\to\mathcal{D}^*[\alpha]$:
$$
P\to\varphi_{\alpha,\,P},\quad
\varphi_{\alpha,\,P}(D)=\xi_{\alpha}(DP),
$$
where $\mathcal{D}^*$ is a graded dual Hopf algebra, and
$P\in\mathcal{P},\;D\in\mathcal{D}$.
\end{defin}
It is easy to see that $\varphi$ is a ring homomorphism. This
follows from the equality
$$
D_{\omega}(P\times Q)=\mu\Delta D_{\omega}(P\otimes Q),
$$
which holds for all compositions $\omega$ and all
polytopes $P,Q$. Here $\mu$ is the multiplication in the ring of polytopes.
\begin{defin}
Let $P^n$ be an $n$-dimensional polytope. Let us define a
quasi-symmetric function
$$
f_m(\alpha,t_1,\dots,t_m)(P^n)=\xi_{\alpha}\Phi(t_m)\dots\Phi(t_1)P^n=\alpha^n+\sum\limits_{k=1}^{\min\{m,n\}}\sum\limits_{0\leqslant
a_1<\dots<a_k\leqslant
n-1}f_{a_1,\,\dots,\,a_k}\alpha^{a_1}M_{(n-a_k,\,\dots,\,a_2-a_1)}.
$$
\end{defin}
For $n=1$ we obtain
$f_1(\alpha,t_1)(P^n)=\alpha^n+\sum\limits_{i=1}^nf_i\alpha^it_1^{n-i}$
is a homogeneous $f$-polynomial in two variables (\cite{Buch}). So the polynomial
$f_n$ is a generalization of the $f$-polynomial.

Consider the increasing sequence of rings
$$
\mathbb Z[\alpha]\subset \mathbb Z[\alpha,t_1]\subset\mathbb
Z[\alpha,t_1,t_2]\subset\dots
$$
with the restriction maps $u_m:\mathbb
Z[\alpha,t_1,\dots,t_{m+1}]\to\mathbb Z[\alpha,t_1,\dots,t_m]$
$$
(u_mg)(\alpha,t_1,\dots,t_m)=g(\alpha,t_1,\dots,t_m,0).
$$
Since
$f_{m+1}(\alpha,t_1,\dots,t_m,0)(P^n)=f_m(\alpha,t_1,\dots,t_m)(P^n)$,
we obtain the ring homomorphism
\begin{gather*}
f:\mathcal{P}\to\Qsym[t_1,t_2,\dots][\alpha]\subset\varprojlim\limits_m\mathbb
Z[\alpha,t_1,\dots,t_m]:\\
f(\alpha,t_1,t_2,\dots)(P^n)=\alpha^n+\sum\limits_{k=1}^{n}\sum\limits_{0\leqslant
a_1<\dots<a_k\leqslant
n-1}f_{a_1,\,\dots,\,a_k}\alpha^{a_1}M_{(n-a_k,\,\dots,\,a_2-a_1)}.
\end{gather*}

It follows from the formula, that the restriction
$$
f(\alpha,t_1,\dots,t_m,t_{m+1},\dots)\to
f(\alpha,t_1,\dots,t_m,0,0,\dots)=f_m(\alpha,t_1,\dots,t_m)
$$
is injective on the space of all $n$-dimensional polytopes,
$n\leqslant m$.
\begin{thm} The image of the space $\mathcal{P}^n$ generated by all
$n$-dimensional polytopes in the ring\linebreak
$\Qsym[t_1,\dots,t_m][\alpha],\;m\geqslant n$ under the mapping $f_m$ consists
of all the homogeneous polynomials of degree $n$  satisfying the equations
\begin{enumerate}
\item
\begin{gather*}
f(\alpha,t_1,-t_1,t_3,\dots,t_m)=f(\alpha,0,0,t_3,\dots,t_m);\\
f(\alpha,t_1,t_2,-t_2,t_4,\dots,t_m)=f(\alpha,t_1,0,0,t_4\dots,t_m);\\
\dots\\
f(\alpha,t_1,\dots,t_{m-2},t_{m-1},-t_{m-1})=f(\alpha,t_1,\dots,t_{m-2},0,0);
\end{gather*}
\item
$$
f(-\alpha,t_1,\dots,t_{m-1},\alpha)=f(\alpha,t_1,\dots,t_{m-1},0);
$$
\end{enumerate}
These equations are equivalent to the Bayer-Billera (generalized
Dehn-Sommerville) relations.\label{f}
\end{thm}
\begin{proof}
$\Phi(-t)\Phi(t)=1=\Phi(0)\Phi(0)$, therefore
\begin{gather*}
\xi_{\alpha}\Phi(t_m)\dots\Phi(t_3)\Phi(-t_1)\Phi(t_1)P^n=\xi_{\alpha}\Phi(t_m)\dots\Phi(t_3)\Phi(0)\Phi(0)P^n;\\
\xi_{\alpha}\Phi(t_m)\dots\Phi(t_4)\Phi(-t_2)\Phi(t_2)\Phi(t_1)P^n=\xi_{\alpha}\Phi(t_m)\dots\Phi(t_4)\Phi(0)\Phi(0)\Phi(t_1)P^n;\\
\dots\\
\xi_{\alpha}\Phi(-t_{m-1})\Phi(t_{m-1})\Phi(t_{n-2})\dots\Phi(t_1)P^n=\xi_{\alpha}\Phi(0)\Phi(0)\Phi(t_{m-2})\dots\Phi(t_1)P^n.
\end{gather*}
On the other hand, Proposition \ref{Euler} gives the last relation
$$
\xi_{-\alpha}\Phi(\alpha)\Phi(t_{m-1})\dots\Phi(t_1)P^n=\xi_{\alpha}\Phi(0)\Phi(t_{m-1})\dots\Phi(t_1)P^n.
$$

Now let us proof the opposite inclusion, that is if the homogeneous polynomial $g$ of degree $n$ satisfies the conditions of the theorem, then $g=f_m(p^n)$ for some $p^n\in\mathcal{P}^n$.

\begin{lemma}
Let $P^n$ be an $n$-dimensional polytope. Then the equation
$$
f_m(\alpha,t_1,\dots,t_q,-t_q,\dots,t_m)=f_m(\alpha,t_1,\dots,0,0,\dots,t_m)
$$
is equivalent to the generalized Dehn-Sommerville relations
$$
\sum\limits_{j=a_t}^{a_{t+1}-1}(-1)^{j-a_t-1}f_{a_1,\,\dots,\,a_t,\,j,\,a_{t+1},\,\dots,\,a_k}=\left(1+(-1)^{a_{t+1}-a_t}\right)f_{a_1,\,\dots,\,a_t,\,a_{t+1},\,\dots,\,a_k}
$$
for $1\leqslant k\leqslant\min(m-1,n-1),\;k+1-q\leqslant t\leqslant m-q$.

The equation
$$
f_m(-\alpha,t_1,\dots,t_{m-1},\alpha)=f_m(\alpha,t_1,\dots,t_{m-1},0)
$$
is equivalent to the generalized Dehn-Sommerville relations
$$
\sum\limits_{j=0}^{a_1-1}f_{j,\,a_1,\,\dots,\,a_k}=(1+(-1)^{a_1-1})f_{a_1,\,\dots,\,a_k}
$$
for $0\leqslant k\leqslant\min(m-1,n-1)$.\label{equations}
\end{lemma}
\begin{proof}
\begin{multline*}
f_m(\alpha,t_1,\dots,t_q,-t_q,\dots,t_m)=\alpha^n+\sum\limits_{k=1}^{\min(m,\,n)}\sum\limits_{0\leqslant
a_1<\dots<a_k\leqslant
n-1}f_{a_1,\,\dots,\,a_k}\alpha^{a_1}\left(\sum\limits_{1\leqslant
l_1<\dots<l_k\leqslant m}t_{l_1}^{n-a_k}\dots
t_{l_k}^{a_2-a_1}\right)=\\
=\alpha^n+\sum\limits_{k=1}^{\min(m-2,\,n)}\sum\limits_{0\leqslant
a_1<\dots<a_k\leqslant
n-1}f_{a_1,\,\dots,\,a_k}\alpha^{a_1}\left(\sum\limits_{1\leqslant
l_1<\dots<l_k\leqslant m,\,l_j\ne q,\,q+1}t_{l_1}^{n-a_k}\dots
t_{l_k}^{a_2-a_1}\right)+\\
+\sum\limits_{k=1}^{\min(m-1,\,n)}\sum\limits_{0\leqslant
a_1<\dots<a_k\leqslant
n-1}f_{a_1,\,\dots,\,a_k}\alpha^{a_1}\left(\sum\limits_{1\leqslant
l_1<\dots<l_j=q<\dots<l_k\leqslant m,\,l_{j+1}\ne
q+1}t_{l_1}^{n-a_k}\dots t_q^{a_{k+2-j}-a_{k+1-j}}\dots
t_{l_k}^{a_2-a_1}\right)+\\
+\sum\limits_{k=1}^{\min(m-1,\,n)}\sum\limits_{0\leqslant
a_1<\dots<a_k\leqslant
n-1}f_{a_1,\,\dots,\,a_k}\alpha^{a_1}\left(\sum\limits_{1\leqslant
l_1<\dots<l_{j+1}=q+1<\dots<l_k\leqslant m,\,l_j\ne
q}t_{l_1}^{n-a_k}\dots (-t_q)^{a_{k+1-j}-a_{k-j}}\dots
t_{l_k}^{a_2-a_1}\right)+\\
\sum\limits_{k=1}^{\min(m,\,n)}\!\!\!\!\!\sum\limits_{0\leqslant
a_1<\dots<a_k\leqslant
n-1}f_{a_1,\,\dots,\,a_k}\alpha^{a_1}\left(\sum\limits_{1\leqslant
l_1<\dots<l_j=q<l_{j+1}=q+1<\dots< l_k\leqslant
m}t_{l_1}^{n-a_k}\dots
(-1)^{a_{k+1-j}-a_{k-j}}t_q^{a_{k+2-j}-a_{k-j}}\dots
t_{l_k}^{a_2-a_1}\right)
\end{multline*}
The first summand is exactly $f(\alpha,t_1,\dots,0,0,\dots,t_m)$. Therefore all the coefficients of the polynomial consisting of the
last three summands should be equal to $0$.

Consider the monomial $\alpha^{a_1}t_{l_1}^{n-a_k}\dots t_q^{a_{t+1}-a_t}\dots t_{l_k}^{a_2-a_1}$. Here $q=l_{k+1-t}$. The existence of a monomial of this form in the sum is equivalent to the conditions
$$
k\leqslant\min(m-1,n-1),\;k+1-t\leqslant q,\;t-1\leqslant m-q-1.
$$
The coefficient of the monomial should be equal to $0$. This is equivalent to the relation:
$$
\left(1+(-1)^{a_{t+1}-a_t}\right)f_{a_1,\,\dots,\,a_t,\,a_{t+1},\,\dots,\,a_k}+\sum\limits_{j=a_t+1}^{a_{t+1}-1}(-1)^{j-a_t}f_{a_1,\,\dots,\,a_t,\,j,\,a_{t+1},\,a_k}=0
$$
Now let us consider the remaining relation
$f(-\alpha,t_1,\dots,t_{m-1},\alpha)=f(\alpha,t_1,\dots,t_{m-1},0)$:
\begin{multline*}
=(-\alpha)^n+\sum\limits_{k=1}^{\min(m-1,\,n)}\sum\limits_{0\leqslant
a_1<\dots<a_k\leqslant
n-1}f_{a_1,\,\dots,\,a_k}(-\alpha)^{a_1}\left(\sum\limits_{1\leqslant
l_1<\dots<l_k\leqslant m-1}t_{l_1}^{n-a_k}\dots
t_{l_k}^{a_2-a_1}\right)+\\
+\sum\limits_{k=1}^{\min(m,\,n)}\sum\limits_{0\leqslant
a_1<\dots<a_k\leqslant
n-1}f_{a_1,\,\dots,\,a_k}(-\alpha)^{a_1}\left(\sum\limits_{1\leqslant
l_1<\dots<l_k=m}t_{l_1}^{n-a_k}\dots
t_{l_{k-1}}^{a_3-a_2}\alpha^{a_2-a_1}\right)=\\
\alpha^n+\sum\limits_{k=1}^{\min(m-1,\,n)}\sum\limits_{0\leqslant
a_1<\dots<a_k\leqslant
n-1}f_{a_1,\,\dots,\,a_k}\alpha^{a_1}\left(\sum\limits_{1\leqslant
l_1<\dots<l_k\leqslant m-1}t_{l_1}^{n-a_k}\dots
t_{l_k}^{a_2-a_1}\right)
\end{multline*}
This is equivalent to the relations
$$
(-1)^{a_1}f_{a_1,\,\dots,\,a_k}+\sum\limits_{j=0}^{a_1-1}(-1)^jf_{j,\,a_1,\,\dots,\,a_k}=f_{a_1,\,\dots,\,a_k}
$$
for $0\leqslant k\leqslant\min(m-1,n-1)$. If $k=0$, then the corresponding relation is exactly the Euler formula
$$
(-1)^n+(-1)^{n-1}f_{n-1}+\dots+f_2-f_1+f_0=1.
$$
\end{proof}
\begin{cor}
For $m\geqslant n$ relations 1. and 2. of the theorem are equivalent to the generalized Dehn-Sommerville relations for the polytope $P^n$:
For $S\subset\{0,\dots,n-1\}$, and $\{i,k\}\subseteq S\cup\{-1,n\}$
such that $i<k-1$ and $S\cap\{i+1,\dots,k-1\}=\varnothing$:
$$
\sum\limits_{j=i+1}^{k-1}(-1)^{j-i-1}f_{S\cup\{j\}}=(1-(-1)^{k-i-1})f_S.
$$
\end{cor}
\begin{proof}
We see that relations 1. and 2. follow from the generalized Dehn-Sommerville relations.

On the other hand, if $i=-1$ then the corresponding relation follows from equation 2.

If $S=\{a_1,\dots,a_t,a_{t+1},\dots,a_s\}$, $i=a_t\geqslant 0$, and $k=a_{t+1}$; or $S=\{a_1,\dots,a_t\}$, $i=a_t$, and $k=n$, then we can take
$q$ such that $s+1-t\leqslant q\leqslant m-t$.
\end{proof}
\begin{rem}
For  $n$-dimensional polytopes for different $m\geqslant n$ not all the equations are independent.
In fact, for $i=-1$ relation (\ref{X2}) follows from equation 2. and all the equations of type 1. do not contain the case $i=-1$.

On the other hand, let $S=\{a_1,\dots,a_t,a_{t+1},\dots,a_s\}$, $i=a_t\geqslant 0$, and $k=a_{t+1}$; or $S=\{a_1,\dots,a_t\}$, $i=a_t$, and $k=n$.
Lemma \ref{equations} implies that relation (\ref{X2}) follows from the equation with $s+1-t\leqslant q\leqslant m-t$. Let us denote $a=s+t-1$, $b=m-t$.
There are two conditions for $s,t$, namely
$$
1\leqslant s\leqslant n-1 \mbox{ and }1\leqslant t\leqslant s.
$$
Let us rewrite these conditions in terms of $a$ and $b$:
$t=m-b;\;s=a+t-1=a-b+m-1$, so
$$
m-n\leqslant b-a\leqslant m-2;\;b\leqslant m-1;\;1\leqslant a
$$
This gives us a triangle on the plane $(a,b)$.
Each point $(a,b)$ of this triangle corresponds to the conditions $a\leqslant q\leqslant b$,
that is there should be the equation for $q$ in the segment $[a,b]$. We can imagine that this
segment is the segment $[(a,a),(a,b)]$ on the plane.

Thus for $m=n$ all the equations for $t_1,\dots,t_{n-1}$ are necessary, for $m=n+1$ it is enough to take the equations for
$q=2,4,6,\dots,2[\frac{n}{2}]$. If $m\geqslant 2n-2$ one equation
$$
f(\alpha,t_1,\dots,t_{n-1},-t_{n-1},t_{n+1},\dots,t_m)=f(\alpha,t_1,\dots,0,0,t_{n+1},\dots,t_m)
$$
gives all the relations of type 1.
\end{rem}
Now let us finish the proof of the theorem.
If the homogeneous polynomial $g\in\Qsym[t_1,\dots,t_m][\alpha]$ of degree $n$ satisfies all the relations of the theorem, then it's coefficients satisfy the generalized Dehn-Sommerville relations. Therefore all the coefficients are linear combinations of the coefficients $g_{a_1,\,\dots,\,a_k}$, where $S=\{a_1,\dots,a_k\}\in \Psi^n$. As we know, the vectors $\{f_{S}(Q),\;S\in\Psi^n\},\;Q\in\Omega^n$ form a basis of the abelian group of all the vectors $\{f_{S},\;S\in\Psi^n,\; f_{S}\in\mathbb Z\}=\mathbb Z^{c_n}$. So the vector $\{g_{S},\;S\in\Psi^n\}$ is an integer combination of the vectors $\{f_{S}(Q),\;S\in\Psi^n\},\; Q\in\Omega^n$. This implies that the polynomial $g$ is an integer combination of the polynomials $f_m(Q),\;Q\in\Omega^n$ with the same coefficients.
\end{proof}

Let us remind that rank of the space $f_m(\mathcal{P}^n),\;m\geqslant n$ is equal to $c_n$.

\begin{prop}
Let $m\geqslant 2$. Then $f_m(\alpha,t_1,\dots,t_m)(P^n)=f_1(\alpha,t_1+\dots+t_m)(P^n)$ if and only if $P^n$ is a simple polytope.
Here $f_1(\alpha,t)$ is a usual homogeneous $f$-polynomial in two variables.
\end{prop}
\begin{proof}
On the ring of simple polytopes
$
\left.d_k\right|_{\mathcal{P}_s}=\left.\frac{d^k}{k!}\right|_{\mathcal{P}_s},
$
so $\Phi(t)=e^{td}$.

Then $\Phi(t_m)\Phi(t_{m-1})\dots\Phi(t_1)=\Phi(t_1+\dots+t_m)$, therefore $f_m(\alpha,t_1,\dots,t_m)(P^n)=f_1(\alpha,t_1+\dots+t_m)(P^n)$.

On the other hand, let $f_m(\alpha,t_1,\dots,t_m)(P^n)=f_1(\alpha,t_1+\dots+t_m)(P^n)$.

Then $f_2(\alpha,t_1,t_2)(P^n)=f_1(\alpha,t_1+t_2)(P_n)$. So
$$
\alpha^n+\sum\limits_{i=0}^{n-1}f_i\alpha^i(t_1^{n-i}+t_2^{n-i})+\sum\limits_{0\leqslant i<j\leqslant n-1}f_{ij}\,\alpha^it_1^{n-j}t_2^{j-i}=\alpha^n+\sum\limits_{i=0}^{n-1}f_i\alpha^i(t_1+t_2)^{n-i}
$$
In particular, $f_{01}=nf_0$, since the coefficients of the monomial $t_1^{n-1}t_2$ on the left and on the right are equal.
Hence $2f_1=f_{01}=nf_0$. This implies that the polytope $P^n$ is simple.
\end{proof}
\begin{rem}
Letting $m$ tend to infinity we obtain that $f(\alpha,t_1,t_2,\dots)(P^n)=f_1(\alpha,t_1+t_2+\dots)(P^n)$ if and only if $P^n$ is simple.
\end{rem}
In the case of simple polytopes the equations of the first type are trivial, but the equation of the
second type has the form:
$$
f_1(-\alpha,t_1+\dots+t_{m-1}+\alpha)=f_1(\alpha,t_1+\dots+t_{m-1})
$$
If we denote $t=t_1+\dots+t_{m-1}$, then
$f_1(-\alpha,\alpha+t)=f_1(\alpha,t)$. This equation is equivalent to
the Dehn-Sommerville relations (after the change of variables
$h(\alpha,t)=f_1(\alpha-t,t)$ it is equivalent to the fact that
$h(\alpha,t)=h(t,\alpha)$)

\section{Characterization of $f$}
In this part we find the condition that uniquely determines the generalized $f$-polynomial.

At first let us find the relation between $f(P^n)$ and $f(d_kP^n)$.
\begin{prop}
For any polytope $P^n\in\mathcal{P}$ we have
$$
f(\alpha,t_1,t_2,\dots)(d_kP^n)=\frac{1}{k!}\left.\frac{\partial^k}{\partial t^k}\right|_{t=0}f(\alpha,t,t_1,t_2,\dots)(P^n)
$$\label{dk}
\end{prop}
\begin{proof}
Indeed,
\begin{multline*}
\frac{1}{k!}\left.\frac{\partial^k}{\partial t^k}\right|_{t=0}f(\alpha,t,t_1,t_2,\dots)(P^n)=\frac{1}{k!}\left.\frac{\partial^k}{\partial t^k}\right|_{t=0}\lim\limits_{r\to\infty}\xi_{\alpha}\Phi(t_r)\dots\Phi(t_1)\Phi(t)P^n=\\
\lim\limits_{r\to\infty}\frac{1}{k!}\left.\frac{\partial^k}{\partial t^k}\right|_{t=0}\xi_{\alpha}\Phi(t_r)\dots\Phi(t_1)\Phi(t)P^n=\lim\limits_{r\to\infty}\xi_{\alpha}\Phi(t_r)\dots\Phi(t_1)d_kP^n=f(\alpha,t_1,t_2,\dots)(d_kP^n)
\end{multline*}
\end{proof}
\begin{cor}
\begin{gather*}
f(\alpha,t,t_1,t_2,\dots)(P^n)=f(\alpha,t_1,t_2,\dots)(P^n)+f(\alpha,t_1,t_2,\dots)(dP^n)t+\\
+f(\alpha,t_1,t_2,\dots)(d_2P^n)t^2+\dots+f(\alpha,t_1,t_2,\dots)(d_nP^n)t^n\\
f_{m+1}(\alpha,t,t_1,\dots,t_m)(P^n)=f_m(\alpha,t_1,\dots,t_m)(P^n)+f_m(\alpha,t_1,\dots,t_m)(dP^n)t+\\
+f_m(\alpha,t_1,\dots,t_m)(d_2P^n)t^2+\dots+f_m(\alpha,t_1,\dots,t_m)(d_nP^n)t^n,\;m\geqslant 0;
\end{gather*}\label{eqn}
\end{cor}
\begin{proof}
The first equality follows from Proposition \ref{dk}.

In fact, both equalities can be proved directly:
\begin{multline*}
f_{m+1}(\alpha,t,t_1,\dots,t_m)(P^n)=\xi_{\alpha}\Phi(t_m)\dots\Phi(t_1)\Phi(t)P^n=\sum\limits_{k=0}^n\left(\xi_{\alpha}\Phi(t_m)\dots\Phi(t_1)d_kP^n\right)t^k=\\
=\sum\limits_{k=0}^{n}f_m(\alpha,t_1,\dots,t_m)(d_kP^n)t^k.
\end{multline*}

Letting $m$ tend to infinity we obtain the first equality.
\end{proof}

The first equality of Corollary \ref{eqn}
is equivalent to the condition $f(\alpha,t_1,t_2,\dots)(\Phi(t)P^n)=f(\alpha,t,t_1,t_2,\dots)(P^n)$.

Thus we see that the following diagram commutes
$$
\begin{CD}
\mathcal{P} @>f>>\Qsym[t_1,t_2,\dots][\alpha]\\
@V{\Phi}VV @VTVV\\
\mathcal{P}[t]@>f>>\Qsym[t,t_1,t_2,\dots][\alpha]
\end{CD}
$$
where $f(t)=t$, and $T:\Qsym[t_1,t_2\dots][\alpha]\to\Qsym[t,t_1,t_2,\dots][\alpha]$ is a ring homomorphism:
$$
Tg(\alpha,t_1,t_2\dots)=g(\alpha,t,t_1,t_2,\dots),\quad g\in\Qsym[t_1,t_2,\dots][\alpha].
$$
\begin{rem}
Note that $f_m(\alpha,t_1,\dots,t_m)(\Phi(t)P^n)=f_{m+1}(\alpha,t,t_1,\dots,t_m)(P^n)$ for all $m\geqslant 0$.

Consider the ring homomorphism
$T_{m+1}:\Qsym[t_1,\dots,t_m][\alpha]\to\Qsym[t,t_1,\dots,t_m][\alpha]$ defined as
$$
T_{m+1}(\alpha)=\alpha,\quad T_{m+1}M_{\omega}(t_1,\dots,t_m)=M_{\omega}(t,t_1,\dots,t_m)
$$
Then the corresponding diagram commutes only for $m\geqslant n$.
\end{rem}

\begin{thm}
Let $\psi:\mathcal{P}\to\Qsym[t_1,t_2,\dots][\alpha]$ be a linear map such that
\begin{enumerate}
\item $\psi(\alpha,0,0,\dots)(P^n)=\alpha^n$;
\item The following diagram commutes:

$$
\begin{CD}
\mathcal{P} @>\psi>>\Qsym[t_1,t_2,\dots][\alpha]\\
@V{\Phi}VV @VTVV\\
\mathcal{P}[t]@>\psi>>\Qsym[t,t_1,t_2,\dots][\alpha]
\end{CD}
$$
\end{enumerate}
Then $\psi=f$.
\end{thm}
\begin{proof}
The first condition implies that $\psi(\alpha,0,0,\dots)(p)=\xi_{\alpha}p$.

Let $\psi(P^n)=\alpha^n+\sum\limits_{\omega}\psi_{\omega}(\alpha)M_{\omega}$, and let $\omega=(j_1,\dots,j_k)$.

Since $\psi(\alpha,t,t_1,t_2,\dots)(P^n)=\psi(\alpha,t_1,t_2,\dots)(\Phi(t)P^n)$, we have
$$
\psi(\alpha,y_1,y_2,\dots,y_k,t_1,t_2,\dots)(P^n)=\psi(\alpha,t_1,t_2,\dots)(\Phi(y_k)\dots\Phi(y_1)P^n)
$$
Therefore
$$
\psi(\alpha,y_1,\dots,y_k,0,0,\dots)=\psi(\alpha,0,0,\dots)(\Phi(y_k)\dots\Phi(y_1)P^n)=\xi_{\alpha}\Phi(y_k)\dots\Phi(y_1)P^n=f(\alpha,y_1,\dots,y_k,0,0,\dots)
$$
Hence $\psi_{\omega}(\alpha)=f_{\omega}\alpha^{n-|\omega|}$ for all $\omega$, so $\psi=f$.
\end{proof}
\begin{rem}
Let us mention that the mapping $T$ is an isomorphism $\Qsym[t_1,t_2,\dots][\alpha]\to\Qsym[t,t_1,t_2,\dots][\alpha]$, while  $\Phi$ is an injection and its image is described by the condition: $p(t)\in\Phi(\mathcal{P})$ is and only if $\Phi(-t)p(t)\in\mathcal{P}$.
\end{rem}
On the ring of simple polytopes we have $f(\alpha,t_1,t_2,\dots)(p)=f_1(\alpha,t_1+t_2+\dots)(p)=f_1(\alpha,\sigma_1)(p)$, so the image of $\mathcal P$ belongs to $\mathbb Z[\alpha,\sigma_1]$. On the other hand, $\Phi(t)=e^{td}$, so
we have the condition $f_1(\alpha,t_1)(e^{td}p)=f(\alpha,t+t_1)(p)$.
\begin{prop}
Let $\psi:\mathcal{P}_s\to\mathbb Z[\alpha,t]$ be a linear mapping such that
\begin{enumerate}
\item $\psi(\alpha,0)(P^n)=\alpha^n$;
\item One of the following equivalent conditions holds:
   \begin{enumerate}
     \item $\psi(\alpha,t_1)(e^{td}p)=\psi(\alpha,t+t_1)(p)$;
     \item $\psi(\alpha,t)(dp)=\frac{\partial}{\partial t}\psi(\alpha,t)(p)$.
    \end{enumerate}
\end{enumerate}
Then $\psi=f_1$.
\end{prop}
\begin{proof}
Let conditions $1$ and $2a$ hold. Then
$$
\psi(\alpha,t+t_1)(p)=\psi(\alpha,t_1)(e^{td}p)=\psi(\alpha,0)(e^{(t+t_1)d}p)=\xi_{\alpha}\Phi(t+t_1)p=f_1(\alpha,t_1+t)(p).
$$

It remains to prove that conditions $2a$ and $2b$ are equivalent.

Indeed, if $\psi(\alpha,t_1)(e^{td}p)=\psi(\alpha,t+t_1)(p)$, then
$\psi(\alpha,t)(p)=\psi(\alpha,0)(e^{td}p)$, therefore
$$
\frac{\partial}{\partial t}\psi(\alpha,t)(p)=\frac{\partial}{\partial t}\psi(\alpha,0)(e^{td}p)=\psi(\alpha,0)(e^{td}dp)=\psi(\alpha,t)(dp).
$$

Now let $\psi(\alpha,t)(dp)=\frac{\partial}{\partial t}\psi(\alpha,t)(p)$.
Then
$$
\frac{\partial}{\partial t}\psi(\alpha,t_1)(e^{td}p)=\psi(\alpha,t_1)(d(e^{td}p))=\frac{\partial}{\partial t_1}\psi(\alpha,t_1)(e^{td}p)
$$
Hence $\left(\frac{\partial}{\partial t}-\frac{\partial}{\partial t_1}\right)\psi(\alpha,t_1)(e^{td}p)=0$, so $\psi(\alpha,t_1)(e^{td}p)$ depends on $\alpha$ and $t+t_1$.

Thus $\psi(\alpha,t_1)(e^{td}p)=\psi(\alpha,t+t_1)(e^{0d}p)=\psi(\alpha,t+t_1)(p)$.
\end{proof}
Proposition 19 in the form $1, 2b$ was first proved in \cite{Buch}.
\section{Ring $\mathcal{D}^*$}
Let us identify $\mathcal{M}=\Qsym[t_1,t_2,\dots]$. Let us remind
that the representation $R:\mathcal{Z}\to\mathcal{D}$ is defined
on the generators as $Z_k\to d_k$. Then $R^*$ is the mapping
$$
R^*:\mathcal{D}^*[\alpha]\to\Qsym[t_1,t_2,\dots][\alpha]:\quad \left(R^*\psi(\alpha)\right)Z^{\omega}=\psi(\alpha)\left(RZ^{\omega}\right).
$$
The ring $\mathcal{D}^*[\alpha]$ is a graded ring with $\deg\alpha^i\psi_i=i+\deg\psi_i$.
\begin{prop}
$$
R^*\varphi_{\alpha,\,P}=f(P).
$$
\end{prop}
For each $m$ let us define the ring homomorphism
$R_m^*:\mathcal{D}^*[\alpha]\to\Qsym[t_1,\dots,t_m][\alpha]$:
$$
R^*_m(\alpha,t_1,\dots,t_m)(\psi(\alpha))=\psi(\alpha)\Phi(t_m)\dots\Phi(t_1)=R^*(\alpha,t_1,\dots,t_m,0,0,\dots).
$$
Then $f_m(p)=R^*_m\varphi_{\alpha,\,p}$ for all $p\in\mathcal{P}$.

The restriction $R^*\to R_m$ is invective for the graduations
$n\leqslant m$ and the mapping $R^*$ is injective on the ring
$\mathcal{D}^*[\alpha]$, since the mapping $R$ is surjective.
\begin{prop}
Let $m\geqslant n$. Then the image of the $n$-th graded component of the
ring $\mathcal{D}^*$ in the ring $\Qsym[t_1,\dots,t_m]$ under the
map $R^*_m$ consists of all the homogeneous polynomials of degree
$n$ satisfying the relations:
\begin{gather*}
g(t_1,-t_1,t_3,\dots,t_m)=g(0,0,t_3,\dots,t_m);\\
g(t_1,t_2,-t_2,t_4,\dots,t_m)=g(t_1,0,0,t_4,\dots,t_m);\\
\dots\\
g(t_1,\dots,t_{m-2},t_{m-1},-t_{m-1})=g(t_1,\dots,t_{m-2},0,0).
\end{gather*}\label{D*}
\end{prop}
\begin{proof}
The ring $\mathcal{D}^*\subset\mathcal Z^*$ consists of all the linear functions $\psi\in\mathcal{Z}^*$ satisfying the property:
$$
\psi z_1\Phi(t)\Phi(-t)z_2=\psi z_1z_2
$$
for all $z_1,z_2\in\mathcal{Z}$, that is the coefficients of all $t^k,\;k\geqslant 1$, on the left are equal to $0$.

Since $R^*_m\psi(\alpha)=\psi(\alpha)\Phi(t_m)\dots\Phi(t_1)$, the relations of the theorem are valid.

On the other hand, let $g\in\Qsym[t_1,\dots,t_m]$ be a homogeneous polynomial of degree $n$ satisfying all the relations of the theorem.
Let us prove that the corresponding linear function $\psi\in\mathcal{M}$ belongs to the image of $\mathcal{D}^*$ under the embedding $R^*$.

It is sufficient to prove the relation (\ref{X8}) in the case when $z_1,z_2$ are monomials.
Since
$$
\Phi(t)\Phi(-t)=1+(Z_1-Z_1)t+(2Z_2-Z_1^2)t^2+\dots=1+(2Z_2-Z_1^2)t^2+\dots,\mbox{ and }\deg\psi=n,
$$
the cases $\deg z_1+\deg z_2=n$ and $n-1$ are trivial.
Let
$$
z_1=Z^{\omega},\;z_2=Z^{\omega'},\;\omega=(j_1,\dots,j_l),\;\omega'=(j_1',\dots,j_{l'}'),\;|\omega|+|\omega'|=n-k,\;k\geqslant 2.$$
Then the only equality we need to prove is
$$
\psi Z^{\omega}\left(\sum\limits_{i=0}^{k}(-1)^iZ_{k-i}Z_i\right)Z^{\omega'}=0
$$
Let us consider the equality
$$
g(t_1,\dots,t_l,t_{l+1},-t_{l+1},t_{l+3},\dots,t_{l+2+l'},\dots,t_m)=g(t_1,\dots,t_l,0,0,t_{l+3},\dots,t_{l+2+l'},\dots,t_m).
$$
The coefficient of the monomial $t_1^{j_1}\dots t_l^{j_l}t_{l+1}^{k}t_{l+3}^{j_1'}\dots t_{l+2+l'}^{j_{l'}'}$ on the left is exactly
$$
\sum\limits_{i=0}^k(-1)^i\psi Z^{\omega}Z_{k-i}Z_iZ^{\omega'}=\psi Z^{\omega}\left(\sum\limits_{i=0}^{k}(-1)^iZ_{k-i}Z_i\right)Z^{\omega'},
$$
and on the right it is equal to $0$.
So $\psi Z^{\omega}\Phi(t)\Phi(-t)Z^{\omega'}=0$ for all $\omega,\;\omega'$, therefore $\psi\in\mathcal{D}^*$.
\end{proof}

\begin{prop}
The image of the space $\mathcal{P}^n$ under the mapping
$\varphi:\mathcal{P}\to\mathcal{D}^*[\alpha]$ consists of all
homogeneous functions $\psi(\alpha)$ of degree $n$ such that
\begin{equation}\label{X8}
\psi(-\alpha)\Phi(\alpha)D=\psi(\alpha)D
\end{equation}
for all $D\in\mathcal{D}$.
\end{prop}
\begin{proof}
According to Proposition \ref{Euler} we have $\xi_{-\alpha}\Phi(\alpha)=\xi_{\alpha}$, so if $\psi(\alpha)=\varphi_{\alpha,\,p}$, then
$$
\psi(-\alpha)\Phi(\alpha)D=\xi_{-\alpha}\Phi(\alpha)Dp=\xi_{\alpha}Dp=\psi(\alpha)D.
$$
On the other hand, let condition (\ref{X8}) hold. Then the polynomial $g(\alpha,t_1,\dots,t_n)=R^*_n\psi(\alpha)$ satisfies the relation
$$
g(-\alpha,t_1,\dots,t_{n-1},\alpha)=\psi(-\alpha)\Phi(\alpha)\Phi(t_{n-1})\dots\Phi(t_1)=\psi(\alpha)\Phi(0)\Phi(t_{n-1})\dots\Phi(t_1)=g(\alpha,t_1,\dots,t_{n-1},0).
$$
Proposition \ref{D*} implies that the polynomial $g$ satisfies all the relations of Theorem \ref{f}, so $g=f_n(p^n)$.
We have $R^*_n\psi(\alpha)=g=f_n(p^n)=R^*_n\varphi_{\alpha,\,p^n}$. Since the restriction $R^*\to R^*_n$ is injective on $j$-th graded component of the ring $\mathcal{D}^*$,  $j\leqslant n$, and $R^*$ is an embedding, we obtain that $\psi(\alpha)=\varphi_{\alpha,\,p^n}$.
\end{proof}
There is a right action of the ring $\mathcal{D}$ on it's graded dual $\mathcal{D}^*$:
$$
(\psi,D)\to\psi D:\quad (\psi D)(D')=\psi(DD').
$$
If $\deg\psi=n$ and $\deg D=k$, then $\deg\psi D=n-k$.
Similarly, there is an action of $\mathcal{D}[[\alpha]]$ on $\mathcal{D}^*[[\alpha]]$.
Then condition (\ref{X8}) means that $\psi(-\alpha)\Phi(\alpha)=\psi(\alpha)$.
Since for any homogeneous function $\psi(\alpha)$ of degree $n$ the function $\psi(\alpha)\Phi(\alpha)$ still has degree $n$, we obtain the corollary.
\begin{cor} Let $\psi(\alpha)\in\mathcal{D}^*[\alpha]$. Then
\begin{equation}\label{X9}
\psi(\alpha)\in \varphi(\mathcal{P})\Leftrightarrow \psi(-\alpha)\Phi(\alpha)=\psi(\alpha).
\end{equation}\label{Im}
\end{cor}
Let $\psi(\alpha)=\psi_0+\psi_1\alpha+\dots+\psi_n\alpha^n$, $u(\alpha)=\psi_0+\psi_2\alpha^2+\psi_4\alpha^4+\dots$, $v=\psi_1\alpha+\psi_3\alpha^3+\dots$. Then $u(\alpha)\Phi(\alpha)-v(\alpha)\Phi(\alpha)=u(\alpha)+v(\alpha)$. Therefore
\begin{equation}\label{X10}
u(\alpha)\left(\Phi(\alpha)-1\right)=v(\alpha)\left(\Phi(\alpha)+1\right).
\end{equation}
For example,
$$
\psi_0d=2\psi_1;\quad \psi_0d_3+\psi_2d=\psi_1d_2+2\psi_3.
$$

Consider the ring $\mathcal{D}^*\otimes\mathbb Z[\frac{1}{2}]$. Then for the graded dual Hopf algebra we have: $(\mathcal{D}\otimes \mathbb Z[\frac{1}{2}])^*\simeq\mathcal{D}^*\otimes\mathbb Z[\frac{1}{2}]$.
The condition that describes the image of the ting $\mathcal{P}\otimes\mathbb  Z[\frac{1}{2}]$ in $\mathcal{D}^*\otimes\mathbb Z[\frac{1}{2}][\alpha]$ is the same,
namely relation (\ref{X10}).

\begin{prop}
In the ring $\mathcal{D}^*\otimes\mathbb Z[\frac{1}{2}][\alpha]$ relation (\ref{X10}) is equivalent to the relation:
\begin{equation}\label{X11}
v(\alpha)=u(\alpha)\frac{\Phi(\alpha)-1}{\Phi(\alpha)+1}=u(\alpha)\frac{\frac{\Phi(\alpha)-1}{2}}{1+\frac{\Phi(\alpha)-1}{2}}=u(\alpha)\sum\limits_{k=1}^{\infty}(-1)^{k-1}\left(\frac{\Phi(\alpha)-1}{2}\right)^k.
\end{equation}\label{uv}
\end{prop}
\begin{proof}
Of course, relation (\ref{X11}) follows from relation (\ref{X10}).
On the other hand, we need to prove, that the function $v(\alpha)$ defined by relation (\ref{X10}) is odd, if $u(\alpha)$ is even.
Indeed,
\begin{gather*}
v(-\alpha)=u(-\alpha)\frac{\Phi(-\alpha)-1}{\Phi(-\alpha)+1}=u(\alpha)\frac{\frac{1}{\Phi(\alpha)}-1}{\frac{1}{\Phi(\alpha)}+1}=u(\alpha)(1-\Phi(\alpha))\Phi(\alpha)^{-1}(1+\Phi(\alpha))^{-1}\Phi(\alpha)=\\
=-u(\alpha)\frac{\Phi(\alpha)-1}{\Phi(\alpha)+1}=-v(\alpha).
\end{gather*}
\end{proof}

Let us take $\varphi_{\alpha,\,p}$ with $\alpha=0$. Then we obtain the classical map
$$
\varphi_0:\mathcal{P}\to\mathcal{D}^*:\quad p\to\varphi_{0,\,p}\quad\varphi_0(p)D=\xi_{0}Dp,\;\forall p\in\mathcal{P},\;D\in\mathcal{D}
$$
\begin{prop}
The mapping $\varphi_0\otimes 1:\mathcal{P}\otimes\mathbb
Z[\frac{1}{2}]\to\mathcal{D}^*\otimes\mathbb Z[\frac{1}{2}]$ is a surjection.
\end{prop}
\begin{proof}
Let $\psi_0\in\mathcal{D}^*\otimes\mathbb Z[\frac{1}{2}]$. Consider the element $\psi(\alpha)=u(\alpha)+v(\alpha)$, $u(\alpha)=\psi_0$, $v(\alpha)=\psi_0\frac{\Phi(\alpha)-1}{\Phi(\alpha)+1}$. According to $\mathbb Z[\frac{1}{2}]$-version of Corollary \ref{Im} we obtain that
$\psi(\alpha)=\varphi_{\alpha,\,p}\in\varphi\otimes 1(\mathcal{P}\otimes\mathbb Z[\frac{1}{2}])$.
Then $\psi_0=\psi(0)=\varphi_{0,\,p}$.
\end{proof}
\begin{quest}
What is the image of the map $\varphi_0$ over the integers?
\end{quest}

\begin{rem}
Let us note, that the space $\mathcal{P}/\Ker{\varphi}$
consists of the equivalence classes of integer combinations of
polytopes under the equivalence relation: $p\sim q$ if and only if
$p$ and $q$ have equal flag f-vectors. Rank of the $n$-th
graded component of the group $\mathcal{P}/\Ker{\varphi}$ is equal to
$c_n$ (\cite{BB}, see also Section 6).

On the other hand, by Corollary \ref{Fib} rank of the $n$-th graded component of the ring
$\mathcal{D}^*$ is equal to $c_{n-1}$ , so the mapping
$\varphi_0:\mathcal{P}/\Ker\varphi\to\mathcal{D}^*$ is not
injective.
\end{rem}
\begin{example}
Let us consider small dimensions.
\begin{itemize}
\item $n=1$.  $\mathcal{P}/\Ker\varphi$ is generated by $CC=I$. The ring $\mathcal{D}$ in this graduation is generated by
$d$.
Then $\varphi_0(I)=2d^*$.
\item $n=2$. $\mathcal{P}/\Ker\varphi$ is generated by $CCC=\Delta^2$ and $BCC=I^2$. $\mathcal{D}$ has one generator $d_2$.
Then $\varphi_0(\Delta^2)=3d_2^*$, $\varphi_0(I^2)=4d_2^*$. Therefore $d_2^*=\varphi_0(I^2-\Delta^2)$, $\Ker\varphi_0$ is generated by $3I^2-4\Delta^2$
\item $n=3$. $\mathcal{P}/\Ker\varphi$ is generated by $BCCC, CBCC=CI^2, CCCC=\Delta^3$, while $\mathcal{D}$ is generated by $d_3,d_2d$.
Then
$$
\varphi_0(BC^3)=5d_3^*+18(d_2d)^*,\quad\varphi_0(CBC^2)=5d_3^*+16(d_2d)^*,\quad \varphi_0(C^4)=4d_3^*+12(d_2d)^*
$$
So $\im \varphi_0$ has the basis $d_3^*,\,2(d_2d)^*$ and $\Ker\varphi_0$ is generated by $2BC^3-6CBC^2+5C^4$.
\end{itemize}
Here by $D_{\omega}^*$ we denote the element of the dual basis: $\langle D_{\omega}^*,D_{\sigma}\rangle=\delta_{\omega,\,\sigma}$
\end{example}
\section{Multiplicative Structure}
\begin{thm}
The ring $f(\mathcal{P})\otimes\mathbb Q$ is a free polynomial
algebra.
\end{thm}
\begin{proof}
According to Proposition \ref{S} and Corollary \ref{Poly} to the shuffle algebra structure theorem we have:
$$
\mathcal{D}^*\otimes\mathbb Q\simeq\mathbb Q[\Lyn_{odd}],
$$
Where $\Lyn_{odd}$ are Lyndon words consisting of odd positive integers.
This is a free polynomial algebra. Therefore $\mathcal{D}^*\otimes\mathbb Q[\alpha]$ is a free polynomial algebra in the generators $\Lyn_{odd}$ and $\alpha$.

Let $\{f_{\lambda}\}$ be the elements of $\mathcal{D}^*\otimes\mathbb Q$ corresponding to the Lyndon words.
Consider the functions $f_{\lambda}(\alpha)$ defined as
$$
f_{\lambda}(\alpha)=u_{\lambda}(\alpha)+v_{\lambda}(\alpha);\quad u_{\lambda}(\alpha)=\psi_{\lambda},\quad v_{\lambda}(\alpha)=u_{\lambda}(\alpha)\frac{\Phi(\alpha)-1}{\Phi(\alpha)+1}=\psi_{\lambda}\frac{e^{s(\alpha)}-1}{e^{s(\alpha)}+1}
$$
Let us identify $\mathcal{D}^*$ with its image in $\Qsym[t_1,t_2,\dots]$ under the embedding $R^*$.

Then $f_{\lambda}(\alpha)\in f(\mathcal{P}\otimes\mathbb Q)=f(\mathcal{P})\otimes\mathbb Q$.
\begin{lemma}
$f(\mathcal{P})\otimes\mathbb Q=\mathbb Q[f_{\lambda},\alpha^2]$.
\end{lemma}
\begin{proof}
At first let us proof that the polynomials $\{f_{\lambda}\},\alpha^2$ are algebraically independent.

Let $g\in\mathbb Q[x,y_1,y_2,\dots,y_s]$ be a polynomial such that $g(\alpha^2,f_1(\alpha),\dots,f_s(\alpha))=0$.
We can write this equality as
$$
g_0(f_1(\alpha),\dots,f_s(\alpha))+g_2(f_1(\alpha),\dots,f_s(\alpha))\alpha^2+\dots+g_{2t}(f_1(\alpha),\dots,f_s(\alpha))\alpha^{2t}=0
$$
for some integers $s,t$.
Let us set $\alpha=0$. Then we obtain
$$
g_0(f_1,\dots,f_s)=0
$$
Since $\{f_{\lambda}\}$ are algebraically independent, $g_0\equiv 0$.
Let us divide the equality by $\alpha^2$ to obtain
$$
g_2(f_1(\alpha),\dots,f_s(\alpha))+g_2(f_1(\alpha),\dots,f_s(\alpha))\alpha^2+\dots+g_{2t}(f_1(\alpha),\dots,f_s(\alpha))\alpha^{2t-2}=0
$$
Iterating this step in the end we obtain that all the polynomials $g_{2i}\equiv 0$, $i=0,1,\dots,2t$, so $g\equiv0$.

Now let us prove that the polynomials $\{f_{\lambda}(\alpha)\},\alpha^2$ generate $f(\mathcal{P})\otimes \mathbb Q$.

Let $\psi(\alpha)=\psi_0+\psi_1\alpha+\dots+\psi_n\alpha^n\in f(\mathcal{P})\otimes\mathbb Q$.
The polynomials $\{f_{\lambda}\}$ generate $\mathcal{D}^*$, so we can find a polynomial $g_0\in\mathbb Q[y_1,\dots,y_s]$ such that
$$
g_0(f_1,\dots,f_s)=\psi_0
$$
Let us consider the element $\psi(\alpha)-g_0(f_1(\alpha),\dots,f_s(\alpha))\in f(\mathcal{P})\otimes \mathbb Q$.
Since $\theta_1=\theta_0\frac{d}{2}$ for any
$$
\theta(\alpha)=\theta_0+\theta_1\alpha+\theta_2\alpha_2+\dots+\theta_r\alpha^r\in f(\mathcal{P})\otimes \mathbb Q,
$$
we obtain that
$$
\psi(\alpha)-g_0(f_1(\alpha),\dots,f_s(\alpha))=\hat\psi(\alpha)\alpha^2.
$$
$\hat\psi(\alpha)\alpha^2\in f(\mathcal{P})\otimes\mathbb Q$, so $\hat\psi(-\alpha)\alpha^2\Phi(\alpha)=\hat\psi(\alpha)\alpha^2$. Since $\mathcal{D}^*\otimes\mathbb Q[\alpha]$ is a polynomial ring, $\hat\psi(-\alpha)\Phi(\alpha)=\hat\psi(\alpha)$, so $\hat\psi(\alpha)\in f(\mathcal{P})\otimes \mathbb Q$.

Iterating this argument in the end we obtain the expression of $\psi(\alpha)$ as a polynomial in $\{f_{\lambda}\}$ and $\alpha^2$.
\end{proof}
This lemma proves the theorem.
\end{proof}

Let us note that dimension of the $n$-th graded component of the
ring $f(\mathcal{P})\otimes\mathbb Q$ is equal to the $n$-th
Fibonacci number $c_n$.

\begin{cor} There is an isomorphism of free polynomial algebras
$$
f(\mathcal{P})\otimes\mathbb Q\simeq \mathbb Q[\Lyn_{odd},\alpha^2]=\mathcal{N}_{odd}\otimes\mathbb Q[\alpha^2]\simeq \mathbb Q[\Lyn_{12}]=\mathcal{N}_{12}\otimes\mathbb Q.
$$\label{three}
\end{cor}
\begin{proof}
Indeed, these three algebras are free polynomial algebras with the same dimensions of the corresponding graded components: dimension of the $n$-th graded component of each algebra is equal to $c_n$
\end{proof}
Let us denote by $k_n$ the number of multiplicative generators of degree $n$. According to Corollary \ref{three} for $n\geqslant 3$ the number of generators $k_n$ is equal to the number of Lyndon words in $\Lyn_{odd}$ or in $\Lyn_{12}$ of degree $n$.

For example, for small $n$ we have
\begin{center}
\begin{tabular}{|c|c|c|}
\hline
n&$\Lyn^n_{12}$&$\Lyn^n_{odd}$\\
\hline
3&[1,2]&[3]\\
4&[1,1,2]&[1,3]\\
5&[1,2,2],\;[1,1,1,2]&[5],\;[1,1,3]\\
6&[1,1,2,2],\;[1,1,1,1,2]&[1,5],\;[1,1,1,3]\\
7&[1,2,2,2],\;[1,1,2,1,2],\;[1,1,1,2,2],\;[1,1,1,1,1,2]&[7],\;[1,1,5],\;[1,3,3],\;[1,1,1,1,3],\;\\
\hline
\end{tabular}
\end{center}

\begin{cor}
\begin{enumerate}
\item $k_{n+1}\geqslant k_{n}$
\item $k_n\geqslant N_n-2$, where $N_n$ is the number of different decompositions of
$n$ into the sum of odd numbers.
\end{enumerate}
\end{cor}
\begin{rem}
It is a well-known fact, that the number of the decompositions of $n$ into the sum of odd numbers is equal to the number of the decompositions of
$n$ into the sum of different numbers. For example,
$$
\begin{matrix}
2&=&1+1&=&2\\
3&=&1+1+1,\;3&=&1+2,\;3\\
4&=&1+1+1+1,\;1+3&=&1+3,\;4\\
5&=&1+1+1+1+1,\;1+1+3,\;5&=&1+4,\;2+3,\;5
\end{matrix}
$$
\end{rem}
\begin{proof}
\begin{enumerate}
\item
$k_1=k_2=k_3=1$.

Let $n\geqslant 3$ and
let $w=[a_1,\dots,a_k]\in\Lyn^n_{12}$. Then $1w=[1,a_1,\dots,a_k]\in\Lyn^{n+1}_{12}$.
Indeed, for any proper tail $[a_i,\dots,a_k],\; i\geqslant 1$ if $a_i>1$, then $[a_i,\dots,a_k]>1w$. If $a_i=1$ then $i\ne k$ and $[a_{i+1},\dots,a_k]>[a_1,\dots,a_k]$, since $w$ is Lyndon. Thus $k_{n+1}\geqslant k_n$
\item
Any decomposition $n=d_1+\dots+d_k$ into the sum of odd numbers $d_1\leqslant d_2\leqslant\dots\leqslant d_k$ gives a Lyndon word $[d_1,\dots,d_k]\in\Lyn_{odd}$, except for the case, when $d_1=\dots=d_k=d\ne n$.

If $d\geqslant 5$, then the word $[1,d-2,1,d,\dots,d]$ is Lyndon.
When $d=1$ or $d=3$ there can be no Lyndon word corresponding to the decomposition $n=1+\dots+1$ or $n=3+\dots+3$, as it happens for $n=6$:
$$
6=1+5=3+3=1+1+1+3=1+1+1+1
$$
We have $k_6=2$ and $N_6=4$. Thus $k_n\geqslant N_n-2$.
\end{enumerate}
\end{proof}
Then we obtain the corollary:
\begin{cor}
For any $t$ such that $|t|<\frac{\sqrt{5}-1}{2}$ there exists a decomposition of
the generating series of Fibonacci numbers as an infinite product:
$$
\frac{1}{1-t-t^2}=1+t+t^2+2t^3+3t^4+\dots=\sum\limits_{n=0}^{\infty}c_nt^n=\prod\limits_{i=1}^{\infty}\frac{1}{(1-t^i)^{k_i}}
$$
Where the infinite product $\lim\limits_{n\to\infty}\sum\limits_{i=1}^n\frac{1}{(1-t^i)^{k_i}}$ converges absolutely.

The numbers $k_n$ satisfy the properties $k_{n+1}\geqslant k_n\geqslant N_n-2$, where $N_n$ is the number of the decompositions of $n$ into the sum of odd summands.
\end{cor}
\begin{proof}
Let $0\leqslant t<\frac{\sqrt{5}-1}{2}$. Then the series $\sum\limits_{n=0}^{\infty}c_nt^n$
converges. This means that $R_N=\sum\limits_{n=N}^{\infty}c_nt^n\xrightarrow[N\to\infty]{} 0$.
But
$$
S_n=\sum\limits_{k=0}^{\infty}c_kt^k-\prod\limits_{i=1}^{n}\frac{1}{(1-t^i)^{k_i}}=a_{n+1}t^{n+1}+a_{n+2}t^{n+2}+\dots
$$
where $0\leqslant a_i\leqslant c_i$. Then $0\leqslant S_n\leqslant R_n$. Therefore $S_n\xrightarrow[n\to\infty]{}0$.

On the other hand, there is a criterion of an absolute convergence of an infinite product:

{\itshape The product $\prod\limits_{n=1}^{\infty}(1+b_n)$ converges absolutely if and only if the series $\sum\limits_{n=1}^{\infty}b_n$ converges absolutely.}

It follows from the fact that $\frac{|\ln(1+b_n)|}{|b_n|}=|\frac{\ln(1+b_n)}{b_n}|\xrightarrow[n\to\infty]{}1$, when $b_n\xrightarrow[n\to\infty]{}0.$

In our case the series $-\sum\limits_{i=1}^{\infty}k_i\ln(1-t^i)$ converges for all $0\leqslant t<\frac{\sqrt{5}-1}{2}$, therefore $\sum\limits_{i=1}^{\infty}k_it^i$ converges for the same values of $t$. Therefore radius of convergence of this series is at least $\frac{\sqrt{5}-1}{2}$. This implies that the product converges absolutely for $|t|<\frac{\sqrt{5}-1}{2}$.
\end{proof}

Here are the examples in small dimensions:
\begin{center}
\begin{tabular}{|c|c|c|c|c|}
\hline n&generators&additive
basis&$c_n$&$k_n$\\
\hline
1&[1]&[1]&1&1\\

2&[2]&[2],\;[1]*[1]&2&1\\

3&[1,2]&[1,2],\;[2]*[1],\;[1]*[1]*[1]&3&1\\

4&[1,1,2]&[1,1,2],\;[1,2]*[1],\;[2]*[2],\;[2]*[1]*[1],\;[1]*[1]*[1]*[1]*[1]&5&1\\
\hline
5&[1,2,2],\;[1,1,1,2]&[1,2,2],\;[1,1,1,2],\;[1,1,2]*[1],\;[1,2]*[2],\;[1,2]*[1]*[1]&8&2\\
&&[2]*[2]*[1],\;[2]*[1]*[1]*[1],\;[1]*[1]*[1]*[1]*[1]&&\\
\hline
6&[1,1,2,2],\;[1,1,1,1,2]&&13&2\\
7&[1,2,2,2],\;[1,1,2,1,2],\;[1,1,1,2,2],\;[1,1,1,1,1,2]&&21&4\\
\hline
\end{tabular}
\end{center}
The corresponding products are:
\begin{gather*}
(1-t)(1-t^2)=1-t-t^2+t^3;\\
(1-t)(1-t^2)(1-t^3)=1-t-t^2+t^4+t^5-t^6;\\
(1-t)(1-t^2)(1-t^3)(1-t^4)=1-t-t^2+2t^5-t^8-t^9+t^{10};\\
(1-t)(1-t^2)(1-t^3)(1-t^4)(1-t^5)^2=1-t-t^2+2t^6+2t^7-t^8-t^9-2t^{10}-\\
-t^{11}-t^{12}+2t^{13}+2t^{14}-t^{18}-t^{19}+t^{20}.
\end{gather*}
The numbers $k_i$ can be found in the following way:
\begin{gather*}
-\log(1-t-t^2)=-\sum\limits_{i=1}^{\infty}k_i\log(1-t^i);\\
\sum\limits_{n=1}^{\infty}\frac{1}{n}\sum\limits_{j=0}^n{n\choose
j}t^{n+j}=\sum\limits_{i=1}^{\infty}k_i\sum\limits_{r=1}^{\infty}\frac{t^{ri}}{r}
\end{gather*}
Thus for any $N$
$$
\sum\limits_{j=0}^{[\frac{N}{2}]}\frac{{N-j\choose
j}}{N-j}=\frac{1}{N}\sum\limits_{i|N}ik_i.
$$
According to the M\"obius inversion formula we obtain
$$
Nk_N=\sum\limits_{d|N}d\sum\limits_{j=0}^{[\frac{d}{2}]}\frac{{d-j\choose
j}}{d-j}\mu\left(\frac{N}{d}\right),
$$
Where $\mu(n)$ is the M\"obius function, that is
$$
\mu(n)=
\begin{cases}
1,&n=1;\\
(-1)^r,&n=p_1\dots p_r,\;\{p_i\}\mbox{ -- distinct prime
numbers};\\
0,&n\mbox{ is not square-free}.
\end{cases}
$$
For example, if $N$ is a prime number $N=p$, then
$$
pk_p=-1+p\sum\limits_{j=0}^{[\frac{p}{2}]}\frac{{p-j\choose
j}}{p-j};\quad
k_p=\sum\limits_{j=1}^{[\frac{p}{2}]}\frac{{p-j\choose j}}{p-j}.
$$
Each summand is integer, since ${p-j\choose j}=\frac{{p-j-1\choose
j-1}(p-j)}{j}$ is integer, and $(p-j)$ and $j$ are relatively
prime numbers.

Thus $k_5=2,\; k_7=4,\;k_{11}=18$, and so on.

\section{Bayer-Billera Ring}
Let us consider the free graded abelian group
$BB\subset\mathcal{P}$, generated by $1$ and all the polytopes
$Q\in\Omega^n,\;n=1,2,\dots$. Rank of the $n$-th graded component of
this group is equal to $c_n$.

Since the determinant if the matrix $K^n$ is equal to $1$, the
generalized $f$-polynomials $\{f(Q),\;Q\in\Omega^n\}$ form a basis
of the $n$-th graded component of the ring $f(\mathcal{P})$.

So the composition of the inclusion $i:BB\subset\mathcal{P}$ and
the mapping $f:\mathcal{P}\to\Qsym[t_1,t_2,\dots][\alpha]$ is an
isomorphism of the abelian groups $BB$ and $f(\mathcal{P})$. This
gives a projection $\pi:\mathcal{P}\to BB$
$$
\pi(p)=x\in BB:\quad f(p)=f(x).
$$
It follows from the definition, that $\pi\circ i=1$ on the space
$BB$.
\begin{thm}
The projection $\pi:\mathcal{P}\to BB$ defined by the relation
$f(\pi p)=f(p)$ gives the group $BB$ the structure of a
commutative associative ring with the multiplication
$x*y=\pi(x\times y)$, such that $f(x*y)=f(x)f(y)$ and
$BB\otimes\mathbb Q$ is a free polynomial algebra in a countable
set of variables.
\end{thm}

\section{Rota-Hopf Algebra}
\subsection{Hopf Algebra of Graded Posets}
\begin{defin}
Let $P$ be a finite poset with a minimal element $\hat 0$ and a maximal element $\hat 1$.
An element $y$ in $P$ {\itshape covers} another element $x$ in $P$, if $x<y$ and there is no $z$ in $P$ such that $x<z<y$.
A poset $P$ is called {\itshape graded}, if there exists a rank function $\rho:\mathcal{P}\to\mathbb Z$ such that
$\rho(\hat 0)=0$ and $\rho(x)+1=\rho(y)$ if $y$ covers $x$. Let us denote $\rho(x,y)=\rho(y)-\rho(x)$ for $x\leqslant y$, and $\rho(P)=\rho(\hat 1)$.
Two finite graded poset are isomorphic if there exists an order preserving bijection between them.

Let $\mathcal{R}$ be the vector space over the field $\mathbbm k$ with basis the set of all isomorphism classes of finite graded posets.
$\mathcal{R}$ is a graded connected Hopf algebra where the degree of $P$ is its rank $\rho(P)$, the multiplication is a cartesian product of posets
$$
P\cdot Q=P\times Q,\quad (x,y)_{P\times Q}\leqslant (u,v)\mbox{ if and only if } x\leqslant u\mbox{ and } y\leqslant v,
$$
the unit element is the poset with one element $\hat 0=\hat 1$, the comultiplication is
$$
\Delta(P)=\sum\limits_{\hat 0\leqslant z\leqslant \hat 1}[\hat 0, z]\otimes[z,\hat 1],
$$
were $[x,y]$ is the subposet $\{z\in P|x\leqslant z\leqslant y\}$, and the counit
$$
\varepsilon(P)=
\begin{cases}
1,&\mbox{ if } \hat 0=\hat 1;\\
0,& \mbox{ else}
\end{cases}
$$
\end{defin}

The antipode of the algebra $\mathcal{R}$ is given by the formula \cite{Sch1}:
$$
S(P)=\sum\limits_{k\geqslant 0}\sum\limits_{\hat 0=x_0<x_1<\dots<x_k=\hat 1}(-1)^k[x_0,x_1]\cdot[x_1,x_2]\dots[x_{k-1},x_k]
$$
\begin{example}
For the simplest Boolean algebra $B_1=\{\hat 0,\hat 1\}$, which is the face lattice of the point $\pt$
\begin{gather*}
\Delta(B_1)=1\otimes B_1+B_1\otimes 1\\
S(B_1)=-B_1
\end{gather*}
\end{example}
The Hopf algebra of graded posets was originated in the work by Joni and Rota \cite{JR}. Variations of this construction were studied in \cite{Ehr,ABS,Sch1,Sch2}. A generalization of this algebra can be found in \cite{RS}.

\subsection{Natural Homomorphism from the Ring of Polytopes to the Rota-Hopf algebra}
There is a natural linear mapping $L:\mathcal{P}\to\mathcal{R}$ of degree $+1$ that sends a polytope $P$ to its face lattice $L(P)$.

It is easy to see that this mapping is injective, but it is not a ring homomorphism, since it doesn't preserve a unit: it sends the point $\pt$ to the Boolean algebra $L(\pt)=B_1$, while the unit of $\mathcal{R}$ is a one-element set $\{\hat0=\hat 1\}$.

\begin{rem}
Let us note that $L(P\times Q)\ne L(P)\times L(Q)$, since the face lattice $L(P\times Q)$ contains the empty face $\varnothing$, which can be considered as  $\varnothing\times\varnothing$, but evidently has no faces of the form  $F\times\varnothing$ or $\varnothing\times G$, where $F$ and $G$ are non-empty faces of $P$ and $Q$ respectively.
\end{rem}

\begin{defin}
There is a natural linear mapping $l=L*d*:\mathcal{P}\to\mathcal{R}$ of degree $0$ that sends the polytope to the sum of the face lattices of its vertex figures (Let us recall that the linear operator $*$ sends the polytope $P$ to its polar $P^*$.)
$$
l(\pt)=\{\hat 0\},\quad l(P)=L\left((dP^*)^*\right)=\sum\limits_{v-\mbox{vertex}}L(P/v).
$$
\end{defin}

\begin{example}
To illustrate the mapping $l$ let us consider the sets of simplicial and simple polytopes:
\begin{itemize}
\item For a simplicial polytope $P^n$ each vertex figure is again a simplicial polytope. Combinatorially we have $\partial(P^n/v)=\link_{\partial P^n}(v)$. Therefore the operation $*d*$ preserves the linear space of all simplicial polytopes.

\item For a simple polytope $P^n$ all of its vertex figures are simplices, therefore $*d*(P^n)=f_0(P^n)\Delta^{n-1}$ and $l(P^n)=f_0L(\Delta^{n-1})=f_0\{\hat 0,\hat 1\}^{n}=f_0B_1^n=B_n$, where $f_0(P^n)$ is the number of vertices of $P^n$, and $\{\hat 0,\hat 1\}^n=B_1^n=B_n$ is a Boolean algebra. For example, $l(\Delta^n)=(n+1)B_1^n$ and $L(\Delta^n)=B_1^{n+1}$.
\end{itemize}
\end{example}

\begin{prop}
$l$ is a homomorphism of graded rings.
\end{prop}
\begin{proof}
We have
$$
l(P)=\sum\limits_{v-\mbox{ vertex }}L(P/v)=\sum\limits_{v}[v,P],
$$
where $L/v$ is a vertex figure, and $[v,P]$ is the interval between the vertex $v$ and the polytope $P$ in the face lattice $L(P)$.
Then
$$
l(P\times Q)=\sum\limits_{v\times w}[v\times w,P\times Q]=\sum\limits_{v\times w}[v,P]\times [w,Q]=\left(\sum\limits_{v}[v,P]\right)\cdot\left(\sum\limits_{w}[w,Q]\right)=l(P)\cdot l(Q)
$$
Here $v,w$ are vertexes of $P$ and $Q$ respectively.
\end{proof}

\begin{defin}
Consider the linear span of all Boolean algebras $B_n=\{0,1\}^n=L(\Delta^{n-1}),\quad n\geqslant 0$
$$
\mathcal{B}=\Ls(1, B_1, B_2,\dots)\subset\mathcal{R}.
$$
We have $B_iB_i=B_{i+j}$, so it is a subring in $\mathcal{R}$. Let us denote $x=B_1$.  Since $\Delta x=1\otimes x+x\otimes 1$,
it is a Hopf subalgebra isomorphic to the Hopf algebra $\mathbbm k[x],\;\Delta x=1\otimes x+x\otimes 1$.
\end{defin}
\begin{prop}
The image of the ring of simple polytopes $\mathcal{P}_s$ under the mapping $l$ is a $\mathbb Z$-subalgebra in $\mathcal{B}$ multiplicatively generated by $2x$ and $x^2$, that is
$$
l(\mathcal{P}_s)=\mathbb Z[x_1,x_2]/(x_1^2-4x_2),\quad \deg x_1=1,\deg x_2=2
$$
where $x_1$ corresponds to $2x=2B_1$ and $x_2$ -- to $x^2=B_2$
\end{prop}
\begin{proof}
Indeed, for a simple polytope $P^n$ we have  $l(P^n)=f_0x^n$. On the other hand, we know that for any simple polytope $2f_1=nf_0$, therefore for odd $n$ the number of vertices $f_0$ should be even. So $l(\mathcal{P}_s)\in\mathbb Z[2x,x^2]$.
On the other hand, $2x=l(I)$ and $x^2=l(I^2-\Delta^2)$.
\end{proof}

In \cite{Ehr} Richard Ehrenborg introduced the $F$-quasi-symmetric function of a graded poset of rank $n$
$$
F(P)=\sum\limits_{\hat 0=x_0<x_1<\dots<x_{k+1}=\hat 1}M_{(\rho(x_0,\,x_1),\,\rho(x_1,\,x_2),\,\dots,\,\rho(x_k,\,x_{k+1}))}=\sum\limits_{0<a_1<\dots<a_k<n}f_{a_1,\,\dots,\,a_k}(P)M_{(a_1,\,a_2-a_1,\,\dots,\,n-a_k)}
$$
where the sum ranges over all chains from $\hat 0$ to $\hat 1$, and $f_{a_1,\,\dots,\,a_k}$ are flag numbers.
This mapping induces a Hopf algebra homomorphism $F:\mathcal{R}\to\Qsym[t_1,t_2,\dots]$.

\begin{prop}
For a polytope $P^n$ the following relation holds:
$$
F(l(P^n))^*=f(0,t_1,t_2,\dots)(P^n)
$$
Thus
$$
*\circ F\circ l=\left.f\right|_{\alpha=0}.
$$\label{Fl}
\end{prop}
\begin{proof}
For the point $\pt$ the relation is trivial $F(l(\pt))^*=1=f(0,t_1,t_2,\dots)(\pt)$.

Since in the poset $L(P/v)=[v,P]$ rank of $v$ is equal to its dimension $0$, $\rho(F)=\dim F$, and we have:
\begin{multline*}
F(l(P^n))^*=\sum\limits_{v}F([v,P^n])^*=\sum\limits_{v\subset F^{a_2}\subset\dots\subset F^{a_k}\subset P^n}M_{(n-a_k,\,\dots,\,a_2)}=\\
=\sum\limits_{0<a_2<\dots<a_k<n}f_{0,\,a_2,\,\dots,\,a_k}\alpha^0 M_{(n-a_k,\,\dots,\,a_2-0)}=f(0,t_1,t_2,\dots)(P^n)
\end{multline*}
\end{proof}

\begin{prop}
We have
$$
L=l\circ(C-\frac{1}{2}[I]),
$$\label{CI}
where $C$ is a cone operator and $[I]$ is the multiplication by the interval $I$.
So the images of the maps $L$ and $l$ over $\mathbb Z[\frac{1}{2}]$ in dimensions $n\geqslant 1$ coincide.\label{LP}
\end{prop}
\begin{proof}
Let us mention that $1=\{\hat 0\}$ does not belong to the image of $L$.

For a polytope $P$ we have $L(CP)=\{\hat 0,\hat 1\}\cdot L(P)=L(\pt)\cdot L(P)$.

On the other hand, $l(I)=2\{\hat 0,\hat 1\}=2L(\pt)$. So $L(\pt)$ belongs to the image of $l$ if and only if $2$ is invertable.

Let $P$ be an $n$-dimensional polytope, $n\geqslant 1$. Then
\begin{multline*}
l(CP)=L(*d*CP)=L(*dCP^*)=L(*(P^*+CdP^*))=L(P)+L(C*d*P)=\\
=L(P)+L(\pt)\cdot l(P)=L(P)+\frac{1}{2}l(I)\cdot l(P)=L(P)+l\left(\frac{1}{2}I\times P\right)
\end{multline*}
So $L(P)=l\left(CP-\frac{1}{2}I\times P\right)$.
\end{proof}
\begin{rem}
The homomorphism $l$ is not injective: we see that on the ring of simple polytopes it remembers only the number of vertices.
However Proposition \ref{CI} shows that $l$ is invective on the image of the operator $C-\frac{1}{2}[I]$.
\end{rem}
\begin{cor}
$$
F(L(P))^*=\left.Cf(P)\right|_{\alpha=0}-\left.\sigma_1f(P)\right|_{\alpha=0},
$$\label{Ff}
where $\sigma_1=\sum\limits_{i=1}^{\infty}t_i=M_{(1)}$.
\end{cor}
\begin{proof}
This formula follows directly from Proposition \ref{LP}.
\end{proof}
Let us denote the linear mapping $\Qsym[t_1,t_2,\dots][\alpha]\xrightarrow{\left.(C-\sigma_1)\right|_{\alpha=0}}\Qsym[t_1,t_2,\dots]$ by $\Lambda$. Then on elementary monomials $\Lambda$ has a very simple form
$$
\Lambda(\alpha^{a_1}M_{(n-a_k,\,\dots,\,a_2-a_1)})=M_{(n-a_k,\dots,\,a_2-a_1,\,a_1+1)}
$$
\begin{example}
For a simple polytope $P^n$ we have $l(P^n)=f_0\cdot L(\Delta^{n-1})=f_0x^n$,
$$
F(l(P^n))=F(l(P^n))^*=f_0\sigma_1^n=f(0,t_1,t_2,\dots)(P^n)=f_1(0,\sigma_1)(P^n).
$$
Using Corollary \ref{Ff} we obtain that
\begin{multline*}
F(L(P^n))^*=\left.(C-\sigma_1)f_1(\alpha,\sigma_1)(P^n)\right|_{\alpha=0}=\sum\limits_{i=1}^{\infty}t_if_1(t_i,t_1+\dots+t_{i-1})(P^n)=\\
=\sum\limits_{j=1}^{\infty}t_j^{n+1}+\sum\limits_{i=1}^{n}f_{n-i}\left(\sum\limits_{j=1}^{\infty}t_j^{n+1-i}(t_1+\dots+t_{j-1})^i\right)=\\
=M_{(n+1)}+\sum\limits_{i=1}^{n}f_{n-i}\sum\limits_{i=a_1<\dots<a_k<n}{i\choose n-a_k,\cdots, a_3-a_2, a_2-i}M_{(n-a_k,\,\dots,\,a_2-i,\,n+1-i)}
\end{multline*}

For example, $F(L(\Delta^n))^*=\sigma_1^{n+1}$.
\end{example}

Thus the relation between the Ehrenborg $F$-quasi-symmetric function and the generalized $f$-polynomial can be illustrated by two commutative diagrams
$$
\begin{CD}
\mathcal{P} @>f>>\Qsym[t_1,t_2,\dots][\alpha]\\
@V{l}VV @V{\alpha=0}VV\\
\mathcal{R}@>F^*>>\Qsym[t_1,t_2,\dots]
\end{CD}
\hspace{3cm}
\begin{CD}
\mathcal{P} @>f>>\Qsym[t_1,t_2,\dots][\alpha]\\
@V{L}VV @V{\Lambda}VV\\
\mathcal{R}@>F^*>>\Qsym[t_1,t_2,\dots]
\end{CD}
$$
For other relations see the next subsection.
\subsection{Hopf Comodule Structure}
\begin{defin}
By a {\it Hopf comodule} (or {\it Milnor comodule}) over a Hopf
algebra $X$ we mean a $k$-algebra $M$ with a unit provided $M$ is a
comodule over $X$ with a coaction $b\colon M\to X\otimes M$ such that
$b(uv)=b(u)b(v)$, i.e. such that $b$ is a homomorphism of rings.
\end{defin}

The ring homomorphism $l$ can be extended to a right graded Hopf comodule structure on $\mathcal{P}$
\begin{prop}
The homomorphism $\Delta:\mathcal{P}\to\mathcal{P}\otimes\mathcal{R}$:
$$
\Delta(P^n)=\sum\limits_{F\subseteq P^n}F\otimes L(P^n/F)=\sum\limits_{F\subseteq P^n}F\otimes[F,P^n],
$$
where by definition $L(P/P)=[P,P]$, defines on $\mathcal{P}$ a right graded Hopf comodule structure over $\mathcal{R}$ such that
$$
(\varepsilon\otimes \id)\Delta P=1\otimes l(P)
$$
\end{prop}
Here by $k$ we mean $\mathbb Z$ or $\mathbb Q$. In the latter case we should take $\mathcal{P}\otimes\mathbb Q$.
\begin{proof}
We should proof that the following three diagrams commute
$$
\xymatrix{
\mathcal{P}\ar[d]^{\Delta}\ar[r]^{\Delta}&\mathcal{P}\otimes\mathcal{R}\ar[d]^{\id\otimes\Delta}\\
\mathcal{P}\otimes\mathcal{R}\ar[r]^{\Delta\otimes\id}&\mathcal{P}\otimes\mathcal{R}\\
}
\hspace{2cm}
\xymatrix{
\mathcal{P}\ar[d]^{\Delta}\ar[dr]^{\id\otimes 1}&\\
\mathcal{P}\otimes \mathcal{R}\ar[r]^{\id\otimes \varepsilon}&\mathcal{P}\otimes\mathbbm k\\
}
\hspace{2cm}
\xymatrix{
\mathcal{P}\ar[d]^{\Delta}\ar[r]^{l}&R\ar[d]^{1\otimes\id}\\
\mathcal{P}\otimes\mathcal{R}\ar[r]^{\varepsilon\otimes\id}&\mathbb Z\otimes\mathcal{R}\\
}
$$
and that $\Delta$ is a ring homomorphism.
Indeed, we have
\begin{multline*}
(\Delta\otimes\id)\Delta P=(\Delta\otimes\id)\left(\sum\limits_{F\subseteq P}F\otimes[F,P]\right)=\sum\limits_{F\subseteq P}\left(\sum\limits_{G\subseteq F}G\otimes[G,F]\right)\otimes[F,P]=\\
=\sum\limits_{G\subseteq P}G\otimes\left(\sum\limits_{G\subseteq F\subseteq P}[G,F]\otimes[F,P]\right)=\sum\limits_{G\subseteq P}G\otimes\Delta[G,P]
=(\id\otimes\Delta)\Delta P
\end{multline*}
$$
(\id\otimes\varepsilon)\Delta(P)=\sum\limits_{F\subseteq P}F\otimes\varepsilon([F,P])=P\otimes 1,
$$
and
$$
(\varepsilon\otimes 1)\Delta P=\sum\limits_{F\subseteq P}\varepsilon(F)\otimes[F,P]=\sum\limits_{v-\mbox{ vertex }}1\otimes[v,P]=1\otimes\left(\sum\limits_{v}[v,P]\right)=1\otimes l(P).
$$
Let us check that $\Delta$ is a homomorphism of rings $\mathcal{P}\to\mathcal{P}\otimes\mathcal{R}$:
\begin{multline*}
\Delta(PQ)=\sum\limits_{F\times G\subseteq P\times Q}F\times G\otimes[F\times G,P\times Q]=
\sum\limits_{F\subseteq P,\,G\subseteq Q}F\times G\otimes[F,P]\cdot[G,Q]=\\
=\left(\sum\limits_{F\subseteq P}F\otimes[F,P]\right)\cdot\left(\sum\limits_{G\subseteq Q}G\otimes[G,Q]\right)=\Delta(P)\cdot\Delta(Q)
\end{multline*}
$$
\Delta(\pt)=\pt\otimes [\pt,\pt]=1\otimes 1
$$
\end{proof}

\begin{cor}
\begin{itemize}
\item Any ring homomorphism $\varphi:\mathcal{P}\to\mathbbm k$ defines a ring homomorphism $\mathcal{P}\to\mathcal{R}$
$$
\mathcal{P}\xrightarrow{\Delta}\mathcal{P}\otimes\mathcal{R}\xrightarrow{\varphi\otimes\id}\mathbbm k\otimes\mathcal{R}\simeq\mathcal{R}
$$
\item Any homomorphism of abelian groups $\psi:\mathcal{R}\to\mathbb Z$ defines a linear operator $\Psi\in \mathcal{L}(\mathcal{P})$
$$
\mathcal{P}\xrightarrow{\Delta}\mathcal{P}\otimes\mathcal{R}\xrightarrow{\id\otimes\psi}\mathcal{P}\otimes\mathbb Z\simeq\mathcal{P}
$$
\item In particular, if $\psi:\mathcal{R}\to\mathbb Z$ is a multiplicative homomorphism, then $\Psi$ is a ring homomorphism.
\end{itemize}
\end{cor}

\begin{example}
Let $\varphi=\xi_{\alpha}$. Then we obtain the ring homomorphism $l_{\alpha}:\mathcal{P}\to\mathcal{R}[\alpha]$ defined as
$$
l_{\alpha}(P^n)=(\xi_{\alpha}\otimes\id)\Delta P^n=\sum\limits_{F\subseteq P}\alpha^{\dim F}[F,P]
$$
\begin{itemize}
\item If we set $\alpha=0$, then we obtain a usual homomorphism $l$.
\item On the ring of simple polytopes $\mathcal{P}_s$
we have
$$
l_{\alpha}(P^n)=\sum\limits_{F\subseteq P^n}\alpha^{\dim F}\{\hat 0,\hat 1\}^{n-\dim F}=\sum_{F\subseteq P^n}\alpha^{\dim F}x^{n-\dim F}=f_1(\alpha, x)
$$
is a homogeneous $f$-polynomial in two variables.
\end{itemize}
\end{example}

Set $F(\alpha)=\alpha$. Then we have the ring homomorphism $F\colon\mathcal{R}[\alpha]\to\Qsym[t_1,t_1,\dots][\alpha]$.
\begin{prop}
Let $P^n$ be an $n$-dimensional polytope. Then
$$
F(l_{\alpha}(P^n))=f(P^n)^*
$$
\end{prop}
\begin{proof}
The proof is similar to the proof of Proposition \ref{Fl}:
\begin{multline*}
F(l_{\alpha}(P^n))=\sum\limits_{F^{a_1}\subseteq P}\alpha^{a_1}\sum\limits_{F^{a_1}\subset F^{a_2}\subset\dots\subset F^{a_k}\subset P^n}M_{(a_2-a_1,\,\dots,\,n-a_k)}=\sum\limits_{F^{a_1}\subset\dots\subset F^{a_k}\subset P^n}\alpha^{a_1}M_{(a_2-a_1,\,\dots,\,n-a_k)}=f(P^n)^*
\end{multline*}
\end{proof}

\begin{prop}
The following diagram commutes:
$$
\begin{CD}
\mathcal{P}@>f^*>>\Qsym[t_1,t_2,\dots][\alpha]\\
@V{\Delta}VV @V{\Delta}VV\\
\mathcal{P}\otimes\mathcal{R}@>f^*\otimes F>>\Qsym[t_1,t_2,\dots][\alpha]\otimes\Qsym[t_1,t_2,\dots]
\end{CD}
$$
\end{prop}
\begin{proof}
We know that
$$
\Delta(M_{(b_1,\,\dots,\,b_k)})=1\otimes M_{(b_1,\,\dots,\,b_k)}+M_{(b_1)}\otimes M_{(b_2,\,\dots,\,b_k)}+\dots+M_{(b_1,\,\dots,\,b_{k-1})}\otimes M_{(b_k)}+M_{(b_1,\,\dots,\,b_k)}\otimes 1
$$
and
$$
f(\alpha,t_1,t_2,\dots)(P^n)^*=\alpha^n+\sum\limits_{F^{a_1}\subset\dots\subset F^{a_k}}\alpha^{a_1}M_{(a_2-a_1,\,\dots,\,n-a_k)}=\sum\limits_{F^{a_1}\subset\dots\subset F^{a_k}\subset P^n}\alpha^{a_1}M_{(a_2-a_1,\,\dots,\,n-a_k)}
$$
Then
\begin{multline*}
\Delta(f(P^n)^*)=\sum\limits_{F^{a_1}\subset\dots\subset G^{a_s}\subset\dots\subset F^{a_k}\subset P^n}\alpha^{a_1}M_{(a_2-a_1,\,\dots,\,a_s-a_{s-1})}\otimes M_{(a_{s+1}-a_s,\,\dots,\,n-a_k)}=\\
=\sum\limits_{G^{a_s}\subseteq P^n}\left(\sum\limits_{F^{a_1}\subset\dots\subset F^{a_{s-1}}\subset G^{a_s}}\alpha^{a_1}M_{(a_2-a_1,\,\dots,\,a_s-a_{s-1})}\right)\otimes\left(\sum\limits_{G^{a_s}\subset F^{a_{s+1}}\subset\dots\subset F^{a_k}\subset P^n}M_{(a_{s+1}-a_s,\,\dots,\,n-a_k)}\right)=\\
=\sum\limits_{G\subseteq P^n}f^*(G)\otimes F([G,P^n])=(f^*\otimes F)\Delta P^n.
\end{multline*}
\end{proof}

These two propositions summarize the relations between the Ehrenborg $F$-quasi-symmetric function and the generalized $f$-polynomial.
\section{Problem of the Description of Flag Vectors of
Polytopes} We have mentioned that on the space of simple polytopes
$$
f(\alpha,t_1,t_2,\dots)(P^n)=f_1(\alpha,t_1+t_2+\dots)(P^n)=f_1(\alpha,\sigma_1)(P^n).
$$
The only linear relation on the polynomial $f_1$ is
$f_1(-\alpha,\alpha+t)=f_1(\alpha,t)$ and it is equivalent to the
Dehn-Sommerville relations. In fact, for the polynomial
$g=g(\alpha,\sigma_1)\in\Qsym[t_1,t_2,\dots][\alpha]$ this
condition is necessary and sufficient to be an image of an integer
combination of simple polytopes.

One of the outstanding results in the polytope theory is the so-called
$g$-theorem, which was formulated as a conjecture by P.~McMullen \cite{Mc1} in 1970 and proved by R.~Stanley \cite{St1} (necessity) and
L.~Billera and C.~Lee \cite{BL} (sufficiency) in 1980.

For an $n$-dimensional polytope $P^n$ let us define an {\itshape
$h$-polynomial}:
$$
h(\alpha,t)(P^n)=h_0\alpha^n+h_1\alpha^{n-1}t+\dots+h_{n-1}\alpha
t^{n-1}+h_nt^n=f_1(\alpha-t,t)(P^n)
$$

Since $f_1(-\alpha,\alpha+t)=f_1(\alpha,t)$, the $h$-polynomial is
symmetric: $h(\alpha,t)=h(t,\alpha)$. So $h_i=h_{n-i}$.

\begin{defin}
A {\itshape $g$-vector} is a set of numbers $g_0=1$,
$g_i=h_i-h_{i-1}, 1\leqslant i\leqslant [\frac{n}{2}]$.
\end{defin}

For any positive integers $a$ and $i$ there exists a unique
''binomial $i$-decomposition'' of $a$:
$$a=\binom{a_i}{i}+\binom{a_{i-1}}{i-1}+\dots+\binom{a_j}{j},$$
where $a_i>a_{i-1}>\dots>a_j\geqslant j\geqslant 1.$

Let us define
$$a^{\langle
i\rangle}=\binom{a_i+1}{i+1}+\binom{a_{i-1}+1}{i}+\dots+\binom{a_j+1}{j+1},
\quad 0^{\langle i\rangle}=0.$$ For example, $a=\binom{a}{1}$, so
$a^{\langle 1\rangle}=\binom{a+1}{2}=\frac{(a+1)a}{2}$; for $i\ge
a$ we have
$a=\binom{i}{i}+\binom{i-1}{i-1}+\dots+\binom{i-a+1}{i-a+1}$, so
$a^{\langle i\rangle}=a$.
\begin{thm*}[$g$-theorem, \cite{St1}, \cite{BL}]
Integer numbers $(g_0,g_1,\dots,g_{[\frac{n}{2}]})$ form a
$g$-vector of some simple $n$-dimensional polytope if and only if
they satisfy the following conditions:
$$
g_0=1,\;0\leqslant g_1,\;0\leqslant g_{i+1}\leqslant g_i^{\langle
i\rangle},\; i=1,2,\dots, \left[\frac{n}{2}\right]-1.
$$
\end{thm*}

As it is mentioned in \cite{Z1} even for $4$-dimensional non-simple polytopes the corresponding problem is extremely hard.

Modern tools used to obtain linear and nonlinear inequalities satisfied by flag $f$-vectors are based on the notion of a $cd$-index \cite{BK,BE1,BE2,EF,ER,L,St3}, a toric $h$-vector \cite{BE,BL,K,St2,Sts} and its generalizations \cite{F1,F2,L},  the ring of chain operators \cite{BH,BLiu,Kalai}.

In \cite{St1} R.~Stanley constructed for a simple polytope $P$ a projective toric variety $X_P$ so that $h_i$ are the even Betti numbers of the (singular) cohomology of $X_P$.
Then the Poincare duality and the Hard Lefschetz  theorem for $X_P$ imply the McMullen conditions for $P$.
 
In \cite{Mc2} P.~McMullen gave a purely geometric proof of these conditions using the notion of the Polytope algebra.

In \cite{St2} Stanley generalized the definition of the $h$-vector to an arbitrary polytope $P$ in such a way that if the polytope is rational the generalized  $\hat h_i$ are the intersection cohomology Betti numbers of the associated toric variety (see the validation of the claim in \cite{Fi}).

The set if numbers $(\hat h_0,\dots,\hat h_n)$ is called a {\itshape toric} $h$-vector.  It is nonnegative and symmetric: $\hat h_i=\hat h_{n-i}$, and in the case of rational polytope the Hard Lefschetz theorem in the intersection cohomology of the associated toric variety proves the unimodality
$$
1=\hat h_0\leqslant \hat h_1\dots\hat h_{[\frac{n}{2}]}
$$
Generalizing the geometrical methods of P.~McMullen \cite{Mc2} Kalle Karu \cite{K} proved the unimodality for an arbitrary polytope.  

The toric $h$-vector consists of linear combinations of the flag $f$-numbers. In general case it does not contain the full information about the flag $f$-vector. In \cite{L} C.~Lee introduced an ''extended toric'' $h$-vector that carries all the information about the flag $f$-numbers and consists of a collection of nonnegative symmetric unimodal vectors. In \cite{F1} and \cite{F2} J.~Fine introduced the generalized $h$- and $g$-vectors that contain toric $h$- and $g$-vectors as subsets, but can have negative entries.

Let us look on the problem of the description of flag $f$-vectors from the point of view of the ring of polytopes.

In the case of simple polytopes we can use the $g$-theorem to describe all the polynomials\linebreak
$\psi\in\Qsym[t_1,t_2,\dots][\alpha]$ that are
images of simple polytopes under the ring homomorphism $f$:
$$
\psi(\alpha,t_1,t_2,\dots)=f(\alpha,t_1,t_2,\dots)(P^n),\quad
P^n\mbox{ -- a simple $n$-dimensional polytope}.
$$

In general case Theorem \ref{f} gives a criterion for the polynomial
$\psi\in\Qsym[t_1,t_2,\dots][\alpha]$ to be an image of an integer
combination of polytopes.
\begin{quest}
Find a criterion for the polynomial
$\psi\in\Qsym[t_1,t_2,\dots][\alpha]$ to be an image of an
$n$-dimensional polytope.
\end{quest}

\end{document}